\documentclass[12pt]{amsart}

\usepackage{amsmath,amssymb,amscd,latexsym}
\usepackage[sans]{dsfont}

\title{Generators of simple Lie algebras \\ in arbitrary characteristics}
\author{Jean-Marie Bois}
\thanks{The author was partially supported by the European Community's Human Potential Programme under contract HPRN-CT-2002-00287 (RTN Network ``K-Theory, Algebraic Groups and Related Structures'') and a long-term research grant from the D.A.A.D}

\address{Fakult\"at f\"ur Mathematik, Universit\"at Bielefeld, Postfach 10 01 31, 33501 Bielefeld, Germany}
\email{jbois@math.uni-bielefeld.de}

\newenvironment{thm}{\noindent \bf Theorem.\it}{ \rm }
\newenvironment{prop}{\noindent \bf Proposition.\it}{ \rm }
\newenvironment{lemma}{\noindent \bf Lemma.\it}{ \rm }
\newenvironment{cor}{\noindent \bf Corollary.\it}{ \rm }

\renewcommand{\proof}{\noindent {\bf Proof.} }
\newcommand{\nr}{n$^o$}
\newcommand{\ds}{\displaystyle}
\newcommand{\nbf}{\noindent \bf }

\newcommand{\FF}{\mathbb{F}}
\newcommand{\NN}{\mathbb{N}}
\newcommand{\ZZ}{\mathbb{Z}}
\newcommand{\CC}{\mathbb{C}}

\newcommand{\gtg}{\mathfrak{g}}
\newcommand{\gth}{\mathfrak{h}}

\newcommand{\gtm}{\mathfrak{m}}
\newcommand{\gtp}{\mathfrak{p}}

\renewcommand{\sl}{\mathfrak{sl}}
\newcommand{\gl}{\mathfrak{gl}}
\newcommand{\psl}{\mathfrak{psl}}
\newcommand{\pgl}{\mathfrak{pgl}}
\renewcommand{\sp}{\mathfrak{sp}}

\newcommand{\cO}{\mathcal{O}}
\newcommand{\cM}{\mathcal{M}}

\newcommand{\cF}{\mathcal{F}}
\newcommand{\cG}{\mathcal{G}}
\newcommand{\cD}{\mathcal{D}}

\newcommand{\nn}{\boldsymbol{\underline{n}}}
\newcommand{\un}{\mathds{1}}

\newcommand{\egv}{\operatorname{egv}}
\newcommand{\End}{\operatorname{End}}
\newcommand{\Aut}{\operatorname{Aut}}
\newcommand{\id}{\operatorname{id}}
\newcommand{\ord}{\operatorname{ord}}

\newcommand{\pr}{\operatorname{pr}}
\newcommand{\rk}{\operatorname{rk}}
\newcommand{\Span}{\operatorname{Span}}
\newcommand{\Der}{\operatorname{Der}}
\newcommand{\ad}{\operatorname{ad}}

\renewcommand{\d}{\partial}
\renewcommand{\epsilon}{\varepsilon}

\begin{document}

\renewcommand{\labelitemi}{\textbullet}
\renewcommand{\labelenumi}{\arabic{enumi}.}

\makeatletter \def\revddots{\mathinner{\mkern1mu\raise\p@ \vbox{\kern7\p@\hbox{.}}\mkern2mu \raise4\p@\hbox{.}\mkern2mu\raise7\p@\hbox{.}\mkern1mu}} \makeatother 

\begin{abstract}
In this paper we study the minimal number of generators for simple Lie algebras in characteristic $0$ or $p > 3$. We show that any such algebra can be generated by $2$ elements. We also examine the 'one and a half generation' property, i.e. when every non-zero element can be completed to a generating pair. We show that classical simple algebras have this property, and that the only simple Cartan type algebras of type $W$ which have this property are the Zassenhaus algebras.
\end{abstract}

\maketitle

\section*{Introduction}

The motivation for this study goes back to some classical questions concerning generating sets of finite simple groups. It is a well-known result that finite simple groups are generated by 2 elements (Dixon's conjecture, see for instance \cite{Aschbacher,Dixon,Steinberg}). In \cite{Guralnick}, Guralnick and Kantor were able to improve this result by showing the following stronger property: any finite simple group $G$ can be generated by ``one and a half elements''; in other words, for all $x \in G \smallsetminus (1)$, there exists $y \in G$ such that the pair $\{ x, y\}$ generates $G$. The proof in both cases involved probabilistic arguments and relied crucially on the Classification theorem for finite simple groups. Recently, this strong generation property was used in \cite{Kunyavskii} to give new characterisations of the solvable radical of finite groups.

The object of this paper is to examine the question of generating sets for finite-dimensional simple Lie algebras (\cite[Questions 2.3 and 2.4]{Kunyavskii}, see also \cite[Problem 5]{Premet}). Thanks to the recently completed Classification theorem for finite-dimensional simple Lie algebras in characteristic $\neq 2,3$ \cite{Premet,Strade} this problem seems to be accessible to explicit methods. 

These questions have arisen before in characteristic 0. In \cite{Kuranashi}, Kuranashi shows that semi-simple Lie algebras are generated by 2 elements; this result was then used to study dense embeddings of free groups into Lie groups. In \cite{Ionescu}, Ionescu shows that simple Lie algebras over the complex numbers are generated by 1,5 elements. However, due to the use of exponentials, his method does not carry over to positive characteristic. In the present article, we will use a geometric approach which is valid in arbitrary characteristic.

Pairs of generators of Lie algebras also arise in another context, namely the study of the commuting variety (see for instance \cite{Levy,Richardson} and references therein). Let $L$ be a Lie algebra. The affine space $L \times L$ contains a family of locally closed subsets:
$$\cG_{d}(L) = \left\{ (x,y) \in L \times L\ |\ \dim\, \langle x,y \rangle = d \right\}.$$
The fact that $L$ is generated by 2 elements simply means that $\cG_{\dim\, L}(L) \neq \emptyset$; in this case, we refer to $\cG_{\dim\, L}(L)$ as the {\it generating variety of $L$}. Observe that the commuting variety $\mathcal{C}_L = \{ (x,y) \in L \times L\ |\ [x,y] = 0\}$ is contained in the closure of $\cG_2(L)$. In some sense, one can think of the varieties $\cG_{\leq n}(L)$ as parts of a hierarchy starting at the commuting variety and ending at the generating variety.

Let us recall briefly the classification result in positive characteristic. Let $\FF$ be an algebraically closed field of characteristic $p > 3$. The first class of simple Lie algebras over $\FF$ consists of analogues of simple Lie algebras in characteristic 0, i.e. Lie algebras of simple algebraic groups. Such Lie algebras are called {\it classical}. The second class is made up of  {\it graded Cartan type Lie algebras}, which are analogues of infinite-dimensional simple Lie algebras of complex vector fields. They can be realised as derivation subalgebras of divided power algebras (see Section \ref{dpaasd}). The third class, the {\it filtered Cartan type Lie algebras}, consists essentially of filtered deformations of graded Cartan Lie algebras. Finally, in characteristic $p = 5$, there is a further exceptional family of Lie algebras, the {\it Melikian Lie algebras}, which have no analogue in characteristic $0$.

The main result of this article is the following: \\

{\bf Theorem A.} \it 
Let $L$ be a simple Lie algebra in characteristic $p \neq 2,3$. Then $L$ is generated by 2 elements.
\rm \\

One should stress that all proofs are constructive and provide explicit descriptions of generating pairs. In most cases, one of the elements can be chosen to be regular semi-simple. The question of 1,5-generation is also solved for $L$ classical or graded of Cartan type $W$ (Witt algebras, see Section \ref{dpaasd} for a definition).

This paper is divided into four parts. In the first part, we define the generating sets $\cG_d(L)$ for an arbitrary Lie algebra and briefly examine their first properties. In the second part we deal with the case of classical simple Lie algebras. We prove (Theorem \ref{GeneratingSetsClassicalCase}): \\

{\bf Theorem B.} \it
Let $\gtg$ be a classical simple Lie algebra in characteristic $p > 3$, or a simple Lie algebra in characteristic 0. Then $\gtg$ is generated by 1,5 elements. 
\rm \\

The idea of the proof is as follows. Given an element $x \in \gtg \smallsetminus \{0\}$, there exists a Cartan subalgebra $\gth \subseteq \gtg$ such that $x$ has a non-zero component in all root spaces relative to the action of $\gth$. Then we show that for generic $y \in \gth$, the elements $x , y$ generate $\gtg$. 

The third and fourth parts are devoted to the non-classical simple Lie algebras. In Section 3, we give an explicit construction of the simple graded algebras of Cartan type $W$; we also establish some properties which will be used in the study of the 1,5-generation problem. Then we mention the remaining types of simple Lie algebras and state the Classification theorem of Block-Wilson-Premet-Strade (Theorem \ref{BWPS}). 

In Section 4, we prove Theorem A for non-classical simple Lie algebras. We deal separately with the graded and the filtered case. Let us outline the proofs. If $L$ is graded, then we know that the homogeneous part of degree 0 is essentially a classical algebra, say $\gtg$. The standard Cartan subalgebra $\gth \subseteq \gtg$ acts diagonally on $L$. Let $y \in \gth$ be a generic element, i.e. such that the action of $\gth$ can be retrieved in terms of the adjoint action $\ad\, y$. Then, by using Theorem B and a careful analysis of the weights of $\gth$ in $L$, we construct an element $x$ such that the weight components of $x$, together with $y$, generate $L$. We finish the proof by showing that all weight components of $x$ belong to the subalgebra generated by $x$ and $y$.

In the filtered case, we proceed similarly. The main difficulty arises from the fact that $L$ might not contain suitable tori. We overcome this issue as follows. Let $L$ be simple filtered of Cartan type. The associated graded algebra $Gr(L)$ is essentially a simple graded Cartan type algebra. Let $\gth \subseteq Gr(L)_0$ be a Cartan subalgebra of the homogeneous subalgebra of degree 0. We show that there exists an action of $\gth$ on $L$ such that $L \simeq Gr(L)$ as $\gth$-modules; furthermore this action can be retrieved in terms of $\ad\, y$ for a suitable element $y \in L$. Then one can proceed as in the graded case.

We conclude the study by dealing with the question of 1,5-generation in type $W$. More precisely, we prove the following (see \ref{specder} for a precise definition of $W(m,\nn)$): \\

{\bf Theorem C.} \it
The Lie algebra $W(m,\nn)$ is generated by 1,5 elements if and only if $m = 1$.
\rm \\

The strategy of proof is as follows. First, we deal with the restricted case, i.e. when $\nn = \un = (1,\ldots,1)$. Studying the action of $\ad\, y$ on $W(m,\un)$ allows us to give an explicit upper bound for the dimension of the derived subalgebra of $L = \FF \langle x,y \rangle$ when $x$ is homogeneous of maximal degree (Corollary \ref{cork}). For general $\nn$, let $x,y \in W(m,\nn)$ with $x$ homogeneous of maximal degree, and $L = \FF \langle x,y \rangle$. Using a suitable embedding of $W(m,\nn)$ into some $W(N,\un)$ in such a way that the homogeneous components of maximal degree are preserved, we use the previously established upper bound for $\dim [L,L]$ to show that $[L,L] \neq W(m,\nn)$ when $m \neq 1$. The fact that $W(1,\nn)$ is generated by 1,5 elements is shown by means of an explicit computation.

\section{Generating varieties of Lie algebras}

In this section, $\FF$ is an algebraically closed field of arbitrary characteristic and $L$ is a finite-dimensional Lie algebra over $\FF$. All topological notions will refer to the Zariski topology.

\subsection{Stratification of $L \times L$}

\subsubsection{} \label{generalitiesLxL}
Let $x_1,\ldots,x_n \in L$ be elements. We denote by $\FF \langle x_1, \ldots, x_n \rangle \subseteq L$ the Lie subalgebra generated by $x_1,\ldots,x_n$. In the sequel we will be mainly concerned with the case $n = 2$. 

Let $L \times L$ be the cartesian product of $L$ with itself. The automorphism group $\Aut(L)$ acts naturally on $L \times L$ by simultaneous conjugation. Furthermore, the group $GL_2$ also acts on $L \times L$ by the formula:
$$ \left( \begin{array}{cc} a & b \\ c & d \end{array} \right) \cdot (x,y) = (a\, x+b\, y, c\, x + d\, y).$$

\subsubsection{} \label{lesstratesGd}
For all $d \in \NN$, we define a subset $\cG_{d} = \cG_d(L) \subseteq L \times L$ by the following formula:
$$\cG_d = \left\{ (x_1,x_2) \in L \times L\ |\ \dim \FF \langle x_1, x_2 \rangle = d \right\}.$$
Clearly, all subsets $\cG_d$ are stable under the actions of $\Aut(L)$ and of $GL_2$. Note that, a priori, a set $\cG_d$ could be empty even though $\cG_{d+1}$ is not. However, it is easy to see that $\cG_d \neq \emptyset$ for $d = 0,1,2$. \\

\subsubsection{} \label{stratif}
\begin{prop}
Denote $\displaystyle \cG_{\leq d} = \bigcup_{d' \leq d} \cG_{d'}$. For all $d \in \NN$, the subset $\cG_{\leq d}$ is closed in $L \times L$. \\
\end{prop}

\proof\ We will prove the following result. Let $Lib_{\, 2}$ be the free Lie algebra in 2 generators; its elements can be considered as Lie polynomials in 2 variables. Then there exist elements $f_1,\ldots,f_N \in Lib_{\, 2}$ satisfying the following property:
$$(\forall x,y \in L):\ \FF \langle x,y \rangle = \Span_\FF \left( f_i(x,y)\, |, i=1,\ldots,N \right).$$
The important feature of this result is that the Lie polynomials $f_1,\ldots,f_N$ are independent of $x$ and $y$ (they actually depend only on $\dim\, L$). In other words, we can use the same Lie words in $x$ and $y$ to linearly span any subalgebra $\FF \langle x,y \rangle \subseteq L$.

First we show how this claim implies the proposition: namely, each set $\cG_{\leq d}$ can be described as:
$$\cG_{\leq d} = \left\{ (x,y) \in L \times L \ |\ \rk_\FF \{ f_1(x,y), \ldots, f_N(x,y) \} \leq d \right\},$$
which is easily seen to be Zariski closed.

Let us now prove the claim. Let $n = \dim(L)$. For any linear endomorphism $u \in \End(L)$ and $x \in L$, the set $\{ x, u(x), \ldots, u^{n-1}(x)\}$ is stable under $u$. If we take $u = \ad\, y$ for some $y \in L$, this is exactly the sub-$\FF\, y$-module generated by $x$, which we will denote by $X_1$.

Now we define inductively subspaces $X_k$, for $k = 2,3,\ldots$ by the following procedure. Assume that the subspace $X_k$ is defined. Let $[x,X_k]$ be the image of $X_k$ under $\ad\, x$; we define $X_{k+1}$ to be the sub-$\FF\, y$-module generated by $[x,X_k]$. Last, let $X^{(k)} = X_1 + \ldots + X_k$.

This defines an ascending chain of subspaces $X^{(1)} \subseteq X^{(2)} \subseteq \ldots$. Because $L$ has dimension $n$, this chain stops at most at $k = n$. By construction, all $X^{(k)}$ are stable under $\ad\, y$. Also, using the equality $X^{(n+1)} = X^{(n)}$, we can see that $X^{(n)}$ is stable under $\ad\, x$.

Let $M = \FF\, y + X^{(n)}$; then $M$ is stable under $\ad\, y$; since $[x,y] \in X_1 \subseteq X^{(n)}$ it is also stable under $\ad\, x$. Furthermore, by construction we have $x,y \in M$. Thus, $M$ is easily seen to be the smallest subspace of $L$ containing $x$ and $y$, and stable under $\ad\, x$ and $\ad\, y$: therefore, $M = \FF \langle x,y \rangle$.

Now all we need to see is that $M$ is linearly spanned by some fixed set of Lie monomials in $x$ and $y$. From the construction, one can see that $M$ is spanned by $x,y$ and the following:
$$[ \underbrace{ y,[y,\ldots,[y}_{\leq n-1\ \mbox{\tiny times}} ,[ x,[ \underbrace{y,\ldots,[y}_{\leq n-1\ \mbox{\tiny times}},[x,[y,\ldots [y ,x ] \ldots ] ],$$
with at most $n$ occurrences of isolated $x$ and at most $n-1$ occurrences of $y$ in-between. These form an independent finite set of Lie monomials as we wanted. 

\subsubsection{}
\begin{cor}
Let $m$ be maximal such that $\cG_m \neq \emptyset$. The subsets $\{ \cG_0,\cG_1,\ldots,\cG_m\}$ define a stratification of $L \times L$. Each stratum is stable under the natural actions of $\Aut(L)$ and $GL_2$ described in \ref{generalitiesLxL}.

Furthermore, the commuting variety $\mathcal{C}_L = \{ (x,y) \in L \times L\ |\ [x,y] = 0\}$ is contained in the closure of $\cG_2$.
\end{cor}

\subsection{Number of generators}

\subsubsection{} \label{genvar}
Let $L$ be a finite-dimensional Lie algebra. We consider the stratification $\{ \cG_{d} \}_{d \geq 0}$ as defined in \ref{lesstratesGd}. Let $m$ be the maximal integer such that $\cG_m \neq \emptyset$: then $\cG_m$ is a non-empty Zariski open subset of $L \times L$.

If $L$ can be generated by 2 elements, we have $m = \dim(L)$. In this case, we will refer to $\cG_m$ as the {\it generating variety of $L$}. We say that $L$ is generated by {\it one and a half elements} if, for all $x \in L \smallsetminus \{ 0 \}$, there exists $y \in L$ such that $L = \FF \langle x,y \rangle$. We can give a geometric interpretation of this property. Namely, $L$ is generated by 1,5 elements if and only if $\pr_1(\cG_m) = L \smallsetminus \{ 0 \}$, where $\pr_1$ is the first projection $L \times L \to L$.
 
\subsubsection{\bf Generators over non-algebraically closed fields} \label{infinitefields}
Let us assume for the moment that $\FF$ is no longer algebraically closed. Assume however that $\FF$ is infinite. Let $L$ be a Lie algebra over $\FF$ and $\overline{L} = L \otimes_{\FF} \overline{\FF}$ the Lie algebra obtained by extension of scalars to the algebraic closure $\overline{\FF}$. Then $L$ is generated by 2 elements (over $\FF$) if and only if $\overline{L}$ is generated by 2 elements (over $\overline{\FF}$).

Indeed, for $n = \dim(L) = \dim( \overline{L} )$, the generating subset $\cG_n(\overline{L}) \subseteq \overline{L} \times \overline{L}$ is a non-empty Zariski open set. Because $\FF$ is infinite, the subset $L \times L \subseteq \overline{L} \times \overline{L}$ is dense, so that $(L \times L) \cap \cG_n(\overline{L}) \neq \emptyset$. Thus, there exist $x,y \in L$ which generate $\overline{L}$ over $\overline{\FF}$. Since $x$ and $y$ are defined over $\FF$, we have $\dim_\FF\, \FF \langle x,y \rangle = \dim_{\overline{\FF}}\, \overline{\FF} \langle x,y \rangle = n$, so that $x,y$ also generate $L$ over $\FF$.

Note that a similar argument applies to Lie algebras which are generated by 1,5 elements.

\section{Generators of classical simple Lie algebras}

Throughout this section, $\FF$ will denote an algebraically closed field of characteristic 0 or $p > 3$.

\subsection{Classical simple Lie algebras} $ $\\

We briefly recall the construction of classical simple Lie algebras in positive characteristics. These are direct analogues of simple Lie algebras over the complex numbers.

\subsubsection{}  Let $\gtg_\CC$ be a simple Lie algebra over $\CC$ (see for instance \cite{Bourbaki}). Let $\mathcal{C}$ be a Chevalley basis of $\gtg_\CC$ and let $\gtg_\ZZ = \ZZ\, \mathcal{C}$ be the $\ZZ$-span of $\mathcal{C}$. Then $\gtg_\ZZ$ is in fact a Lie algebra over $\ZZ$; we define a Lie algebra over $\FF$ by setting $\gtg_\FF = \gtg_\ZZ \otimes_\ZZ \FF$. 

\subsubsection{}
\begin{prop} In the previous construction, the Lie algebra $\gtg_\FF$ is simple, except if the root system of $\gtg_\CC$ has type $A_n$ and $p \mid n+1$. In this case, $\gtg_\FF \simeq \sl_{n+1}$ has a 1-dimensional centre, and the factor Lie algebra $\gtg_\FF / \mathfrak{z}(\gtg_\FF) \simeq \psl_{n+1}$ is simple. \\
\end{prop}

The simple Lie algebras obtained by this process are called {\it classical}. As in characteristic 0, they are parametrised by Dynkin diagrams of types $A_n$, $B_n$, $C_n$, $D_n$, $G_2$, $F_4$, $E_6$, $E_7$, $E_8$. Types $B,C,D$ correspond to Lie algebras $\mathfrak{so}_{2n+1}$, $\mathfrak{sp}_{2n}$ and $\mathfrak{so}_{2n}$ respectively; type $A$ corresponds to Lie algebras $\sl_{n+1}$ or $\psl_{n+1}$ as in the proposition. Note that, in the modular theory, the term ``classical'' includes the types $E,F,G$.

\subsubsection{}  Consider a classical simple Lie algebra over $\FF$. Let $\gth \subseteq \gtg$ be a Cartan subalgebra. We have a root space decomposition:
\begin{equation} \label{RootSpace}
\gtg = \gth \oplus \bigoplus_{\alpha \in \Phi} \gtg_\alpha,
\end{equation}
where $\Phi \subseteq \gth^*$ is the root system of $\gtg$ associated with $\gth$. This decomposition has the following useful and well-known properties: \\

\begin{lemma} \label{usefulandeasy} 
We retain the above notations. Then:
\begin{enumerate}
\item For all $\alpha \in \Phi$, the subspace $\gtg_{\alpha}$ is of dimension 1.
\item We have $\ds \gth = \sum_{\alpha \in \Phi} [\gtg_{\alpha}, \gtg_{-\alpha}]$. \\
\end{enumerate}
\end{lemma}

\subsection{Generating systems of classical simple Lie algebras} $ $\\

Now we will tackle our first problem, namely the question of generators in classical simple Lie algebras. We show that such algebras have the 1,5-generation property (Theorem \ref{GeneratingSetsClassicalCase}). Corollaries \ref{corproduct} and \ref{corvaria} deal with analogous questions for direct products and central extensions of classical simple Lie algebras. The reason for doing this is that it will apply to the homogeneous part of degree 0 in graded Cartan type Lie algebras (Theorem \ref{GeneratingSetsGradedCase}).

\subsubsection{} 
\begin{lemma} \label{zariskiopeninh}
Let $V$ be a vector space and $\cF \subseteq V^*$ a finite set of linear forms. There exists a non-empty Zariski open subset $\Omega \subseteq V$ such that for all $v \in \Omega$, the elements $\left\{ f(v)\, |\, f \in \cF \right \} \subseteq \FF$ are pairwise distinct. \\
\end{lemma}

In particular, the lemma will be applied in the situation when $V = \gth$ is the Cartan subalgebra of some classical simple Lie algebra and $\cF = \Phi$ is the corresponding root system. \\

\proof\ Write out $\cF = \{ f_1, \ldots, f_n \}$. We identify the symmetric algebra $S(V^*)$ with the algebra of regular maps on $V$. Consider the element $\ds P = \prod_{1 \leq i < j \leq n}(f_j - f_i) \in S(V^*)$. Since all factors $f_j - f_i$ are non-zero, we have $P \neq 0$. Thus, the Zariski open subset $\Omega = \left\{ v \in V\, |\, P(v) \neq 0 \right \} \subseteq V$ is nonempty. By construction, for all $v \in \Omega$ the scalars $\left\{ f(v)\, |\ f \in \cF \right\}$ are pairwise distinct.

\subsubsection{}
\begin{lemma} \label{denseorbit}
Let $\gtg$ be a classical simple Lie algebra over $\FF$, $x \in \gtg$ a non-zero vector and $f_1,\ldots,f_n \in \gtg^*$ be non-zero linear forms on $\gtg$. There exists an automorphism $\sigma \in Aut(\gtg)$ such that $f_i(\sigma(x)) \neq 0$ for all $i \in \{1,\ldots,n\}$. \\
\end{lemma}

\proof\ Let $G$ be a simple connected algebraic $\FF$-group with the same Dynkin diagram as $\gtg$. Then $G$ acts as automorphisms of $\gtg$. Let us check first that $\gtg$ has no non-trivial $G$-stable subspace. Indeed, if $\gtg \not \simeq \mathfrak{psl}_{m}$ with $p \mid m$, then $\gtg \simeq Lie(G)$ and the property follows from the simplicity of $\gtg$. Otherwise, $Lie(G) = \mathfrak{sl}_{m}$, and the centre $\FF\, I_{m}$ is the only nontrivial $G$-invariant subspace. Hence, $\mathfrak{psl}_{m} = \mathfrak{sl}_{m} / \FF\, I_{m}$ has the desired property.

Let us first check the lemma for $n = 1$. For simplicity we write $f = f_1$. Let $X = \Span(G \cdot x) \subseteq \gtg$ be the vector space spanned by the orbit of $x$ under action of $G$. Then $X$ is a non-zero $G$-stable subspace of $\gtg$, hence $X = \gtg$. In particular, the orbit $G \cdot x$ is not contained in the kernel of $f$: in other words, there exists some $g \in G$ such that $f(g.x) \neq 0$.

Now we prove the general case. For all $i \in \{1,\ldots,n\}$, let $U_i(x) = \{g \in G\, |\, f_i(g.x) \neq 0\}$. This is a Zariski open subset of $G$; the previous discussion shows that it is not empty. The group $G$ being connected, the intersection $\ds U(x) = \bigcap_{i=1}^n U_i(x) \subseteq G$ is again a nonempty open subset; by construction, any element $g \in U(x)$ satistfies the required property. The lemma is proved.

\subsubsection{} Now we can address the question of generators in simple classical Lie algebras.\\

\begin{thm} \label{GeneratingSetsClassicalCase}
Let $\FF$ be an algebraically closed field of characteristic $0$ or $p > 3$ and $\gtg$ be a classical simple Lie algebra over $\FF$. Then $\gtg$ is generated by one and a half elements. \\
\end{thm}

\proof\ Let $\gth \subseteq \gtg$ be a Cartan subalgebra and $\Phi \subseteq \gth^*$ be the associated root system. Recall that $\gtg$ has a decomposition $\gtg = \gth \oplus \bigoplus_\alpha \gtg_\alpha$ into root spaces. For $x \in \gtg$, we have a decomposition:
\begin{equation*}
x = x_0 + \sum_{\alpha \in \Phi} x_\alpha,
\end{equation*}
where $x_0 \in \gth$ and each $x_\alpha \in \gtg_\alpha$. 

First we show that applying a suitable automorphism of $\gtg$ allows us to reduce to the case where $x_\alpha \neq 0$ for all $\alpha \in \Phi$. We can choose a basis $\{e_1,\ldots,e_n\}$ of $\gtg$ which is compatible with the root space decomposition. By Lemma \ref{denseorbit} applied to the dual basis $\{e_1^*,\ldots,e_n^*\}$, there exists an automorphism $g \in \Aut(\gtg)$ such that $e_i^*(g.x) \neq 0$ for all $i$. Since each subspace $\gtg_\alpha$ is generated by some element $e_i$, this implies $(g.x)_\alpha \neq 0$ for all $\alpha$ as desired. In the sequel we can assume that all $x_\alpha \neq 0$.

According to Lemma \ref{zariskiopeninh}, there exists a non-empty Zariski open subset $\Omega \subseteq \gth$ such that for all $y \in \Omega$, the elements $\left\{ \alpha(y)\, | \ \alpha \in \Phi \right\}$ are pairwise distinct. Choose any $y \in \Omega$: we will show that $x$ and $y$ generate $\gtg$. Let $L = \FF \langle x,y \rangle \subseteq \gtg$ be the Lie subalgebra generated by $x$ and $y$. For all $k \in \NN$, $L$ contains the element $(\ad\, y)^k (x)$. Hence, we get:
\begin{equation} \label{VdM}
(\ad\, y)^k (x) = \sum_{\alpha \in \Phi} \alpha(y)^k x_\alpha \in L.
\end{equation}
By construction of $y$, the scalars $\{\alpha(y)\, |\, \alpha \in \Phi\}$ are pairwise distinct. So the equations (\ref{VdM}) constitute a Van der Monde invertible system; thus, all components $x_\alpha \in L$. Since $x_\alpha \neq 0$ and $\dim(\gtg_\alpha)=1$ for all $\alpha$, we get $\gtg_\alpha \subseteq L$ for all root $\alpha \in \Phi$. Using Lemma \ref{usefulandeasy}, \nr 2, we also get $\gth \subseteq L$, whence $L = \gtg$: the theorem is proved. \\

\subsubsection{\bf Remark} Actually, one could adapt the proof to gain some insight into the generating variety of simple classical Lie algebras. Recall from \ref{genvar} the generating variety $\cG \subseteq \gtg \times \gtg$. Denote by $\pr_1: \cG \to \gtg$ the restriction to $\cG$ of the first projection map. Also, for $x \in \gtg$, we denote $\cG(x) = \{ y \in \gtg \ |\ \FF\langle x,y \rangle = \gtg\}$, so that the fibre $\pr^{-1}(x) = \{ x \} \times \cG(x)$. Let $\gth$ be a Cartan subalgebra of $\gtg$. If $\gth \cap \cG(x) \neq \emptyset$, then $\cG(x) \supseteq \gth^{reg}$, the regular elements of $\gth$.

\subsubsection{} \label{corproduct}
\begin{cor}
Let $\FF$ be an algebraically closed field of characteristic $0$ or $p > 3$ and $\gtg = \gtg_1 \oplus \ldots \oplus \gtg_n$ be a direct sum of classical simple Lie algebra over $\FF$. Denote by $\pi_k \colon \gtg \twoheadrightarrow \gtg_k$ the natural map, for all $k$. For any element $x \in \gtg$ such that $\pi_k(x) \neq 0$ for all $k$, there exists $y \in \gtg$ such that $\gtg = \FF \langle x , y \rangle$. \\
\end{cor}

\proof\ As before, $\gtg$ has a decomposition $\gtg = \gth \oplus \sum \gtg_\alpha$, where $\gth$ is the sum of Cartan subalgebras of the factors $\gtg_1$, \ldots, $\gtg_n$ and each root space $\gtg_\alpha$ is contained in one of the simple factors $\gtg_k$. Denote by $x = x_0 + \sum x_\alpha$ the decomposition of $x$ with respect to this direct sum. By applying a suitable automorphism of $\gtg$, and using the fact that each $\pi_k(x) \neq 0$, we can reduce to the case where all $x_\alpha \neq 0$. Then we can conclude as before, choosing $y \in \gth$ to be such that the scalars $\alpha(y)$ are pairwise distinct.

\subsubsection{} \label{corvaria}
\begin{cor}
Let $\FF$ be an algebraically closed field of characteristic $0$ or $p > 3$.
\begin{enumerate}
\item Let $\gtg = \gl_n$. For any non-central element $x \in \gtg$, there exists $y \in \gtg$ such that $x$ and $y$ generate $\gtg$.
\item The Lie algebra $\gtg = \mathfrak{pgl}_n$ is generated by 1,5 elements.
\item Let $\gtp$ be a perfect Lie algebra which is generated by 1,5 elements. Let $\gtg$ be a central extension of $\gtp$ with 1-dimensional centre:
\begin{equation} \label{centralextension}
0 \longrightarrow \FF\, z \longrightarrow \gtg \longrightarrow \gtp \longrightarrow 0.
\end{equation} 
For any non-central element $x \in \gtg$, there exists $y \in \gtg$ such that $x$ and $y$ generate~$\gtg$. \\
\end{enumerate}
\end{cor}

\proof\ For part 1, we can use an argument similar to the one used in Theorem \ref{GeneratingSetsClassicalCase}. Statement \nr 2 is a direct consequence of \nr 1.

Let us examine \nr 3. Let $\pi : \gtg \twoheadrightarrow \gtp$ be the canonical projection. Since $x$ is non-central we have $\pi(x) \neq 0$, so that there exists $\eta \in \gtp$ such that $\langle \pi(x), \eta \rangle = \gtp$. Fix an element $y_0 \in \gtg$ such that $\pi(y_0) = \eta$; for all $\alpha \in \FF$ set $y_\alpha = y_0 + \alpha\, z$. Let $L_\alpha = \langle x,y_\alpha \rangle \subseteq \gtg$. We have $\pi(L_\alpha) = \langle \pi(x),\pi(y_\alpha) \rangle = \gtp$. In particular, $L_\alpha$ has codimension at most 1 in $\gtg$.

If $\dim(L_\alpha) = \dim(\gtg)$ for some values of $\alpha$, we are done. Now assume $\dim(L_\alpha) = \dim(\gtg) -1 = \dim(\gtp)$ for all $\alpha$. Then $\pi$ induces an isomorphism $L_\alpha \stackrel{\sim}{\to} \gtp$; hence, the exact sequence (\ref{centralextension}) is split. Thus, we obtain $\gtg = \FF\, z \oplus L_\alpha \simeq \FF\, z \oplus \gtp$. Since $\gtp$ is perfect, we can see that the subalgebra $L_\alpha = [\gtg,\gtg]$ does not depend on $\alpha$. But then $z = y_1 - y_0 \in L_1 + L_0 = L_0$, contradicting $L_\alpha \cap \FF\, z = 0$.

\subsubsection{}
We give some examples of the situations described in Corollary \ref{corvaria}. \\

{\nbf Example 1.} Any direct sum $\gtg = \gtp \oplus \FF\, z$, where $\gtp$ is classical simple and $z$ is central, satisfies the hypotheses of \nr 3. This is the case of the homogeneous subalgebra of degree 0 in the graded Cartan type Lie algebra of type $K$. \\

{\nbf Example 2.} Consider the following diagram:
$$
\begin{CD}
\gl_n @>>>  \sl_n \\ 
@VVV   @VVV \\ 
\pgl_n @>>> \psl_n
\end{CD}
$$
In this diagram, horizontal arrows correspond to taking the derived subalgebra and vertical arrows correspond to factoring out the centre. When $p \mid n$, these four Lie algebras are pairwise non-isomorphic. In this case, statement \nr 1 cannot be retrieved from \nr 3, because $\pgl_n$ is generated by 1,5 elements but is not a perfect Lie algebra, and $\sl_n$ is perfect but not generated by 1,5 elements (because the centre is non-trivial).

\section{Simple Lie algebras in positive characteristics} \label{dpaasd}

For the remainder of the paper, $\FF$ will denote an algebraically closed field of positive characteristic $p > 3$. In this case there exist further classes of finite-dimensional simple Lie algebras, which are analogues of Lie algebras of vector fields over the complex numbers \cite[Section 4.1]{SF}.

In this section we will give a full description of the simple graded Lie algebras of type $W$ as derivations of divided power algebras (Section \ref{GradedTypeW}). We will also mention the other graded types $S, H, K$ (as well as the Melikian algebras in characteristic 5) and outline the construction of simple filtered algebras of Cartan type; however, for technical reasons we will omit the detailed constructions. The interested reader may find these in \cite{Strade}; see also \cite{Premet} for an overview of the classification theorem.

\subsection{Multi-indices and conventions}

\subsubsection{} Let $m \geq 1$ be an integer. A multi-index is an $m$-tuple of integers $\alpha = (\alpha_1,\ldots,\alpha_m)$. We define the factorial and binomial coefficients by the usual conventions: for $\alpha,\beta \in \NN^m$, 
$$\alpha! = \alpha_1 ! \cdots \alpha_m ! \quad \mbox{ and } \quad \binom{\alpha}{\beta} = \binom{\alpha_1}{\beta_1} \cdots \binom{\alpha_m}{\beta_m}.$$
The length of a multi-index $\alpha \in \NN^m$ is defined to be $|\alpha| = \alpha_1+\ldots+\alpha_m$.

Also, define the following distinguished multi-indices: for all $k \in \{1,\ldots,m\}$, the $k$-th unit multi-index is $\epsilon_k = (\delta_{1,k},\ldots,\delta_{m,k})$; we also set $\un = \epsilon_1 + \ldots + \epsilon_m = (1,\ldots,1)$.

\subsubsection{} \label{ordmultind}
We can partially order the set of multi-indices the following way: for $\alpha,\beta \in \ZZ^m$, let $\alpha \leq \beta$ if the components $\alpha_i \leq \beta_i$ for all $i \in \{1,\ldots,m \}$. We will also write $\alpha < \beta$ if $\alpha \leq \beta$ and $\alpha \neq \beta$. 

\subsubsection{\bf Conventions} For a multi-index denoted $\nn \in \NN^m$, we will always make the implicit assumption that $\un \leq \nn$, in other words all the entries of $\nn$ are positive integers. Given such a multi-index $\nn$ we denote by $\tau = \tau(\nn) \in \NN^m$ the multi-index $\tau = (p^{n_1}-1,\ldots,p^{n_m}-1)$.

\subsection{Divided power algebras}

\subsubsection{} Let $m \geq 1$ be an integer. Define the divided power algebra $\cO(m)$ over $\FF$ as follows: as vector space, $\cO(m)$ has a basis $\{ x^{(\alpha)}\, |\, \alpha \in \NN^m \}$ and the multiplication map is given by the formula:
$$(\forall\, \alpha,\beta \in \NN^m)\ :\ x^{(\alpha)} x^{(\beta)} = \binom{\alpha + \beta}{\alpha} x^{(\alpha+\beta)}.$$
It is convenient to set $x^{(\alpha)} = 0$ if one of the components $\alpha_i < 0$. The algebra $\cO(m)$ has a natural grading over $\ZZ$ such that all elements $x^{(\alpha)}$ are homogeneous, of degree $| \alpha |$. \\

{\noindent \bf Remark.} If $\FF$ were of characteristic 0, then $\cO(m)$ would be naturally isomorphic to a polynomial algebra in $m$ variables via the mapping $\cO(m) \stackrel{\sim}{\to} \FF[X_1,\ldots,X_m]$ defined by $x^{(\alpha)} \mapsto \frac{1}{\alpha!} X_1^{\alpha_1} \ldots X_m^{\alpha_m}$. In positive characteristics however, such a map is not available; actually one can show that the algebra $\cO(m)$ is not even finitely generated.

\subsubsection{} \label{structuredp}
Let $\nn \in \NN^m$ be a multi-index with non-zero components. Define a linear subspace:
$$\cO(m,\nn) = \Span \left\{x^{(\alpha)}\ |\ 0 \leq \alpha_i < p^{n_i}\ \mbox{for all}\ i \right\} \subseteq \cO(m).$$

\begin{lemma}
\begin{enumerate}
\item The subspace $\cO(m,\nn)$ is a graded subalgebra of $\cO(m)$, of dimension $p^{|\nn|}$.
\item If $\nn = \un$, the algebra $\cO(m,\nn)$ is isomorphic to the truncated polynomial algebra $B_m = \FF[X_1,\ldots,X_m] / (X_1^p,\ldots,X_m^p)$.
\item For arbitrary $\nn \in \NN^m$, we have $\cO(m,\nn) \simeq \cO(|\nn|,\un)$. \\
\end{enumerate}
\end{lemma}

Items 1 and 2 are stated in \cite[p. 8]{Premet}. For item \nr 3, one could use Harper's theorem on the structure of differentiably simple rings \cite{Harper}; alternatively, one can check that the assignment:
$$(\forall\, i \in \{ 1,\ldots,m\}),\ (\forall j \in \{0,\ldots,n_i-1 \})\ : \ y_{i,j} \mapsto x^{(p^j \epsilon_i)}$$
extends uniquely to an isomorphism $\FF[\{ y_{i,j} \}_{i,j} ] / (y_{i,j}^p)\stackrel{\sim}{\to} \cO(m,\nn)$. 

\subsubsection{\bf Convention.} \label{lanotationBm}
We will often identify the algebras $\cO(m,\un)$ and $B_m$ defined in the lemma. In $B_m$, we will denote by $x_1,\ldots,x_m$ the natural generators, so that $x_i^p = 0$ for all $i$. We will also use the natural convention $x^\alpha = x_1^{\alpha_1} \cdots x_m^{\alpha_m}$. In this case, the multiplication is simply given by the formula $x^\alpha x^\beta = x^{\alpha + \beta}$.

\subsubsection{} \label{usefulofomn}
We can summarise some useful properties of algebras $\cO(m,\nn)$. \\

\begin{prop} $ $
\begin{enumerate}
\item The ring $\cO(m,\nn)$ is a local ring, with maximal ideal $\gtm$ generated by homogeneous elements of non-zero degree.
\item For all $x \in \cO(m,\nn)$, one has $x^p \in \FF$. Furthermore, $x^p = 0$ if and only if $x \in \gtm$.
\item The ring $\cO(m,\nn)$ has a unique minimal non-zero ideal. It is generated by $x^{(\tau)}$, the unique homogeneous element of maximal degree (up to a non-zero scalar).
\item Assume $\nn = \un$. Let $\xi_1,\ldots,\xi_k \in \cO(m,\un)$ be homogeneous elements of degree 1. The following are equivalent: \\
\begin{tabular}{cl} 
{\bf (i)} & $\xi_1,\ldots,\xi_k$ are linearly dependent over $\FF$; \\
{\bf (ii)} & $\xi_1^{p-1} \ldots \xi_k^{p-1} = 0$.
\end{tabular} \\
\end{enumerate}
\end{prop}

\proof\ Using Lemma \ref{structuredp}, we can check all properties for the rings $B_m$. Statements \nr 1, 2 and 3 are easy to check. Let us prove statement \nr 4. First, we check $(i)\ \Rightarrow (ii)$: assume for instance that $\ds \xi_k = \sum_{j=1}^{k-1} \lambda_j\, \xi_j$ with all $\lambda_j \in \FF$. Then:
$$\xi_1^{p-1}\cdots \xi_k^{p-1} = \sum_{j=1}^{k-1} \lambda_j \, \xi_1^{p-1}\cdots \xi_j^p \cdots \xi_{k-1}^{p-1}.$$
Since $\xi_j^p = 0$, we also get $\xi_1^{p-1} \ldots \xi_k^{p-1} =0$ as required.

Conversely, assume that $(i)$ does not hold, i.e. $\xi_1,\ldots, \xi_k$ are linearly independent: we will show that $(ii)$ does not hold. Denote by $V$ the homogeneous subspace of degree 1 in $B_m$ and by $\cD = \sum_{j=1}^m \FF\, D_j \subseteq \Der(B_m)$ the subspace generated by the partial derivatives: then we have a natural identification $\cD \simeq V^*$ (the basis $\{D_1,\ldots,D_m\}$ being dual to $\{x_1,\ldots,x_m\}$).  Thus, there exist $\d_1,\ldots,\d_k \in \cD$ such that $\d_i(\xi_j) = \delta_{i,j}$ for all $i$ and $j$. Then: 
$$\d_1^{p-1}\cdots \d_k^{p-1}(\xi_1^{p-1} \cdots \xi_k^{p-1}) = \d_1^{p-1}(\xi_1^{p-1}) \cdots \d_k^{p-1}(\xi_k^{p-1}) = (p-1)! \cdots (p-1)! \neq 0,$$
so that necessarily $\xi_1^{p-1} \cdots \xi_k^{p-1} \neq 0$.

\subsection{Simple graded Lie algebras of type W} \label{GradedTypeW}

\subsubsection{\bf Special derivations of divided power algebras} \label{specder}
A derivation $D \in \Der\, \cO(m)$ is {\it special} if:
$$(\forall\, \alpha \in \NN^m)\ :\ D(x^{(\alpha)}) = \sum_{i = 1}^m x^{(\alpha-\epsilon_i)} D(x^{(\epsilon_i)}).$$

{\noindent \bf Example 1.} For all $k \in \{1,\ldots,m\}$, define the {\it partial derivative} $D_k$ by the assignment $D_k(x^{(\alpha)}) = x^{(\alpha - \epsilon_k)}$. Then $D_k$ is a special derivation of $\cO(m)$. \\

{\noindent \bf Example 2.} Let $k \in \{1,\ldots,m\}$. Because $\FF$ has characteristic $p > 0$, the $p$-th power $D_k^p$ is again a derivation of $\cO(m)$. It is readily seen that $D_k^p \neq 0$ but $D_k^p(x^{(\epsilon_i)}) = 0$ for all $i \in \{1,\ldots,m\}$: thus $D_k^p$ is not a special derivation of $\cO(m)$.

\subsubsection{} \label{typeW}
Let $W(m) \subseteq \Der\, \cO(m)$ be the Lie algebra of special derivations of $\cO(m)$. It is actually a free $\cO(m)$-module, with basis $\{D_1,\ldots,D_m\}$. We define $W(m,\nn)$ to be the free sub-$\cO(m,\nn)$-module generated by the partial derivatives. Since all partial derivatives stabilise $\cO(m,\nn)$, it is a Lie subalgebra, and we have a natural inclusion:
$$W(m,\nn) \subseteq \Der\, \cO(m,\nn).$$
Actually, it can be checked that the inclusion is proper unless $\nn = \un$. The Lie algebras $W(m,\nn)$ are called {\it graded Lie algebra of Cartan type $W$}, or also {\it Witt algebras}. The following is well-known \cite[Section 4.2]{Strade}: \\

\begin{prop}
For all integer $m \geq 1$ and all $\nn \in \NN^m$, the Lie algebra $W(m,\nn)$ is simple.
\end{prop}

\subsubsection{\bf Natural grading on $W(m,\nn)$} Recall that the associative algebra $\cO(m,\nn)$ has a natural grading for which each element $x^{(\alpha)}$ is homogeneous, of degree $|\alpha|$. Thus, $W(m,\nn)$ inherits a natural grading, for which all derivations of the form $x^{(\alpha)} D_j$ have degree $|\alpha| - 1$. The homogeneous subspace of maximal degree has dimension $m$ and is generated by derivations $x^{(\tau)}\, D_j$, for $\tau$ as in \ref{ordmultind} and $j \in \{1,\ldots,m\}$. The homogeneous subspace of minimal degree -1 has also dimension $m$ and is generated by the partial derivatives $D_j$, for $j \in \{1,\ldots,m\}$.

The gradings on $\cO(m,\nn)$ and $W(m,\nn)$ are compatible, in the sense that if $\delta_i$ is a homogeneous derivation of degree $i$ and $f_j$ a homogeneous element of degree $j$ in $\cO(m,\nn)$, then both $\delta_i(f_j) \in \cO(m,\nn)$ and $f_j\, \delta_i \in W(m,\nn)$ are homogeneous, of degree $i+j$.

\subsubsection{} \label{usefulofwm1}
The following property is an analogue of statement \nr 4 of Proposition \ref{usefulofomn}. As the proof is very similar, it will be omitted. Since $W(m,\un) \subseteq \End\, \cO(m,\un)$, associative products of derivations exist as linear endomorphisms of $\cO(m,\un)$. Recall that $x^{(\tau)}$ is the unique monomial of maximal degree in $\cO(m,\un)$. \\

\begin{lemma} 
Let $\d_1,\ldots,\d_k \in W(m,\un)_{-1}$ be homogeneous derivations of degree $-1$. The following are equivalent: \\
\begin{tabular}{cl} 
{\bf (i)} & $\d_1,\ldots,\d_k$ are linearly dependent over $\FF$; \\
{\bf (ii)} & $\d_1^{p-1} \cdots \d_k^{p-1} = 0$; \\
{\bf (iii)} & $\d_1^{p-1} \cdots \d_k^{p-1} (x^\tau) = 0$.
\end{tabular} \\
\end{lemma}

\subsubsection{} \label{plongement}
According to Lemma \ref{structuredp}, we have $\cO(m,\nn) \simeq \cO(|\nn|,\un)$. This isomorphism induces a natural isomorphism $\Der\, \cO(m,\nn) = \Der\, \cO(|\nn|,\un) = W(|\nn|,\un)$. Thus, we get a natural embedding:
$$\iota: W(m,\nn) \hookrightarrow W(|\nn|,\un).$$
We will usually identify $W(m,\nn)$ with its image, thus considering that $W(m,\nn) \subseteq W(|\nn|, \un)$; however it may be useful sometimes to keep track of the embedding $\iota$ to avoid confusion. Also, one should note that this embedding does not preserve gradings. However, the following proposition shows that it is compatible with the highest degree components: \\

\begin{prop} 
Let $m \geq 1$ be an integer and $\nn \in \NN^m$ a multi-index. Consider the natural gradings of $W(m,\nn)$ and $W(|\nn|,\un)$:
$$W(m,\nn) = \bigoplus_{k=-1}^s W(m,\nn)_k\ \mbox{ and }\ W(|\nn|,\un) = \bigoplus_{k=-1}^{s'} W(|\nn|,\un)_k.$$
Let $\iota: W(m,\nn) \hookrightarrow W(|\nn|,\un)$ be the natural embedding. Then $\iota \big(W(m,\nn)_s\big) \subseteq W(|\nn|,\un)_{s'}$. \\
\end{prop}

\proof\ Recall that $W(m,\nn)_s$ is generated by derivations of the form $x^{(\tau)} D_j$, for $j\in \{1,\ldots,m\}$. Let $J_0 = J_0(m,\nn)$ be the unique minimal ideal of $\cO(m,\nn)$ (see Proposition \ref{usefulofomn}). One can check easily that:
$$W(m,\nn)_s = J_0(m,\nn) \cdot W(m,\nn).$$
We can identify $\cO(m,\nn)$ with $\cO(|\nn|,\un)$. By uniqueness of the minimal ideal $J_0$, we get $J_0(m,\nn) = J_0(|\nn|,\un)$. Thus,
$$W(m,\nn)_s = J_0(m,\nn)\cdot W(m,\nn) \subseteq J_0(|\nn|,\un)\cdot W(|\nn|,\un) =  W(|\nn|,\un)_{s'}.$$
The proposition is proved.

\subsection{Classification of simple Lie algebras in positive characteristics} $ $\\

In this section we briefly mention the other types of simple Lie algebras in positive characteristics. Explicit constructions can be found in \cite{Premet,Strade}; specific properties which are relevant to our study will be recalled when needed (mainly in Section \ref{section3}). We will conclude the section by recalling the classification theorem of Block-Wilson-Premet-Strade (Theorem \ref{BWPS})

\subsubsection{\bf Graded Lie algebras of Cartan type \cite[Section 4.2]{Strade}} 
There exist three other families of simple graded Lie algebras of Cartan type, the {\it special}, {\it hamiltonian} and {\it contact} Lie algebras. These are referred to as simple graded Lie algebras of Cartan type $S,H$ and $K$ respectively. As in type $W$, they can be indexed by a pair $(m,\nn)$ where $m \geq 1$ is an integer and $\nn \in \NN^m$ a multi-index with non-zero entries. Furthermore, there are the following restrictions for $m$:
\begin{itemize}
\item in type $S$, one has $m \geq 3$;
\item in type $H$, the integer $m = 2r$ is even;
\item in type $K$, the integer $m = 2r+1$ is odd.
\end{itemize}

Let $X \in \{ W,S,H,K \}$ be a symbol. The simple Lie algebra of type $X$ attached to $(m,\nn)$ is denoted by $X(m,\nn)^{(\infty)}$; it actually arises as the last term in the derived series of some Lie algebra $X(m,\nn)$.

\subsubsection{\bf Filtered Lie algebras of Cartan type \cite[Section 4.2]{Strade}} 
Another family of simple Lie algebras in characteristic $p > 3$ are the {\it simple filtered Lie algebras of Cartan type}. These are essentially filtered deformations of simple graded Lie algebras of Cartan type.

Let $L$ be a simple filtered Lie algebra of Cartan type. It is endowed with a standard filtration $L = L_{(-s')} \supseteq \ldots \supseteq L_{(s)}$. Let $\cG = Gr(L)$ be the associated graded algebra. Then $\cG$ is not simple in general, but $\cG^{(\infty)}$, the last term in the derived series, is a simple ideal which is isomorphic to a simple graded Lie algebra of Cartan type. By the Compatibility Property \cite[Property (6.2.2)]{Strade}, there exists a symbol $X \in \{ W,S,H,K\}$ such that:
$$X(m,\nn)^{(\infty)} \subseteq Gr(L) \subseteq X(m,\nn);$$
in that case we say that $L$ is of type $X$. The following results are \cite[Corollaries 6.1.7 and 6.6.3]{Strade}: \\

\begin{prop}
Let $L$ be a simple filtered Lie algebra of Cartan type.  If $L$ is of type $W$ or $K$, then $L$ is isomorphic to a graded Lie algebra of same type.
\end{prop}

\subsubsection{\bf Melikian Lie algebras \cite[Section 4.3]{Strade}} 
If $\FF$ has characteristic $p = 5$, there exists another family of graded Lie algebras indexed by pairs of positive integers $\nn = (n_1,n_2)$, the {\it Melikian Lie algebras}. These have no analogue in characteristic 0 or $p > 5$.

\subsubsection{\bf Classification Theorem} \label{BWPS}
The following result is quoted from \cite[Theorem 7]{Premet}: \\

\begin{thm} {\rm (Block-Wilson-Premet-Strade)}

Let $L$ be a finite dimensional simple Lie algebra over an algebraically closed field of characteristic $p > 3$. Then $L$ is either a classical Lie algebra or a filtered Lie algebra of Cartan type or one of the Melikian algebras. 
\end{thm}

\section{Generators of non-classical simple Lie algebras}

Again, in this section, $\FF$ denotes an algebraically closed field of characteristic $p > 3$.

\subsection{Generators of simple graded Lie algebras of Cartan type} \label{section3} $ $ \\

In this section we examine the structure of graded Cartan type algebras and show that they can be generated by two elements. We first recall from \cite{Strade} some relevant properties of these algebras. Section \ref{zerobit} to \ref{torus} are devoted to some useful general features of simple graded Lie algebras of Cartan type. The main source for the results is \cite[Section 4.2]{Strade}. Similar results hold for the Melikian algebras in characteristic 5 \cite[Section 4.3]{Strade}; however to ease exposition we won't mention them explicitly here.

\subsubsection{\bf Notation} \label{zerobit}
In the sequel, we will fix $X \in \{W,S,H,K\}$ a symbol, $m \in \NN$ and $\nn \in \NN^m$ a multi-index. Also, assume that $m \geq 3$ if $X = S$, $m = 2r$ is even if $X = H$ and $m = 2r+1$ is odd if $X = K$. We consider the simple graded Cartan type Lie algebra $L = X(m,\nn)^{(\infty)}$. The structure of the homogeneous subalgebra $L_0$ can be described explicitly \cite[Section 4.2]{Strade}: indeed, $L_0$ is isomorphic to $\gl_m$, $\sl_m$, $\sp_{2r}$, $\sp_{2r} \oplus \FF\, I_{2r}$ for  $X = W,S,H,K$ respectively. In particular, we can apply the constructions of Theorem \ref{GeneratingSetsClassicalCase} to $L_0$. Also, the subalgebras $L_0$ have a natural $p$-mapping arising from associative $p$-th powers of matrices in $\gl_m$.

\subsubsection{} \label{rappels}
Let $L = X(m,\nn)^{(\infty)}$ be a simple graded Lie algebra of Cartan type. Consider its decomposition into homogeneous subspaces:
$$L = L_{-r} \oplus \ldots \oplus L_s,$$
where $r,s \geq 0$. Note that $L_0$ is a subalgebra of $L$, and each $L_i$ is an $L_0$-submodule of $L$. The following can be found in \cite[Section 4.2 and 4.3]{Strade}: \\

\begin{prop}
We keep the general notation above.
\begin{enumerate}
\item The integers $r,s \geq 1$. We have $r = 1$ in type $W,S,H$ and $r = 2$ in type $K$.
\item The subspaces $L_{-r}$, $L_{-1}$ and $L_s$ are irreducible $L_0$-modules.
\item For all $i \in \{-r,\ldots,s-1\}$, one has $L_i = [L_{-1},L_{i+1}]$. In particular, $L$ is generated by $L_{-1}$ and $L_s$.
\end{enumerate}
\end{prop}

\subsubsection{} \label{restrictedaction}
Let $L = \bigoplus_{k \in \ZZ} L_k$ be any graded Lie algebra. Recall that the degree map given by $\deg(x_k) = k\, x_k$ for all homogeneous elements $x_k \in L_k$ extends to a derivation of the Lie algebra $L$. \\

\begin{lemma} Let $L$ be a simple graded Cartan type Lie algebra.
\begin{enumerate}
\item Any homogeneous derivation $\delta$ which vanishes on $L_0$ is proportional to the degree derivation.
\item The subalgebra $L_0$ is restricted, and the adjoint action of $L_0$ on $L$ is a restricted action. In other words, denote by $x^{[p]}$ the image of $x \in L_0$ under the natural $p$-mapping. The relation $(\ad\, x)^p = \ad(x^{[p]})$ holds as derivations of $L$. \\
\end{enumerate}
\end{lemma}

\proof\ Item \nr 1 is a special case of \cite[Proposition 2.1]{Farnsteiner86}. Let us check \nr 2. From \ref{zerobit} we can see that $L_0$ is a restricted Lie algebra. For all $x \in L_0$, the map $(\ad\, x)^p - \ad(x^{[p]}) : L \to L$ is a derivation which vanishes on $L_0$. Using \nr 1 and the fact that $\FF$ is algebraically closed, we see that there exists a map $\chi : L \to \FF$ such that:
$$(\ad\, x)^p - \ad(x^{[p]}) = \chi(x)^p\, \deg.$$
Thus, the adjoint representation of $L_0$ on each $L_k$ has $p$-character $k. \chi$. Therefore, to show that the adjoint representation of $L_0$ in $L$ is restricted it is enough to show that the adjoint representation of $L_0$ on $L_{-1}$ is restricted. This can be shown using a case-by-case examination and the explicit descriptions of these actions given in \cite[Section 4.2]{Strade}.

\subsubsection{} \label{torus}
Following \cite[Section 5.5]{Strade}, there exists a subalgebra $T_X \subseteq L = X(m,\nn)^{(\infty)}$ acting diagonally on $L$. Furthermore, one has $T_X \subseteq L_0$ (it is actually the standard Cartan subalgebra); therefore, $T_X$ also acts diagonally on each homogeneous component of $L$. Denote $\Gamma_k = \Gamma(T_X,L_k) \subseteq T_X^*$ the set of roots of $T_X$ acting on $L_k$. One has:
$$L_k = \bigoplus_{\alpha \in \Gamma_k} L_k^\alpha,$$
and each $L_k^\alpha \neq 0$. We shall need information on the roots of $T_X$ on some of the homogeneous subspaces. \\

\begin{lemma}
With the above notations, we have the following properties:
\begin{enumerate}
\item $\Gamma_0 \cap \Gamma_{-1} = \Gamma_0 \cap \Gamma_s = \emptyset$, where $s$ is the maximal degree of homogeneous components of $L$;
\item there exist $\alpha \in \Gamma_{-1}$, $\beta \in \Gamma_s$ such that $\alpha \neq \beta$. \\
\end{enumerate}
\end{lemma}

\proof\ From \cite[Section 4]{SF}, we can see that the representation of the classical algebra $L_0$ in $L_{-1}$ is contragredient to the natural representation of $L_0$, i.e. the representation given by the embedding $L_0 \hookrightarrow \gl_m$ if $X = W,S,H$ and $L_0 \hookrightarrow \gl_{2r}$ if $X = K$. In each case, the weight spaces are one-dimensional, so that $|\Gamma_{-1}| \geq 2$ whenever $m \neq 1$. In this case, property \nr 2 is obvious. The case $m = 1$ can only occur if $L = W(1,n)$: then, $T_W$ is one-dimensional and acts by the scalar $k$ on each homogeneous component $L_k$. Since $s = p^n-2$ for $W(1,n)$, we get $\Gamma_{-1} = \{-1\}$ and $\Gamma_s = \{-2\}$, yielding the result.

Let us turn to \nr 1. The fact that $\Gamma_0 \cap \Gamma_{-1} = \emptyset$ is \cite[Lemma 5.5.1]{Strade}. To show that $\Gamma_0 \cap \Gamma_{s} = \emptyset$ we use similar arguments to those used in \cite{Strade}, detailed below. \\

{\bf In type $K$:} the centre of $L_0$ acts by the degree. Let $\| \tau \| = \sum_{i=1}^{2r}(p^{n_i} - 1) + 2(p^{n_m}-1)-2$. By \cite[p. 190]{Strade}, the homogeneous component of highest degree for type $K$ has degree $s = \| \tau\| $ if $m+3 \not \equiv 0 \mod p$ and $s = \| \tau \| - 1$ if $m+3 \equiv 0 \mod p$. Since $\| \tau \| \equiv -(m+3) \mod p$, it shows that $s \not \equiv 0 \mod p$ in all cases. Thus, the centre of $L_0$ does not act trivially on $L_s$, forcing $\Gamma_0 \cap \Gamma_s = \emptyset$. \\

For $X \in \{ W,S,H \}$, we will use more explicit descriptions of the involved roots. If $\{h_1,\ldots,h_r\}$ is a basis of $T_X$ (for instance the standard basis of $T_X$ viewed as a Cartan subalgebra of the almost classical Lie algebra $L_0$), denote $\{ \widehat{\epsilon}_1, \ldots, \widehat{\epsilon}_r \} \subseteq T_X^*$ the dual basis. These elements generate an abelian subgroup isomorphic to $(\ZZ / p\, \ZZ)^r$, so that we can always consider that $(\ZZ / p\, \ZZ)^r \subseteq T_X^*$. Also, denote by $\{ \epsilon_1 , \ldots , \epsilon_m \}$ the natural basis of $\ZZ^m$. The following computations can be carried out using the formulas of \cite[Section 4.2]{Strade} or \cite[Chapter 4]{SF}. \\

{\bf In type $W$:} we have $r = m$. Let $\pi_W \colon \ZZ^m \to (\ZZ / p\, \ZZ)^m$ be the natural map. Then we have $\pi_W(\epsilon_i) = \widehat{\epsilon}_i$ for all $i$. Let $\widehat{\tau} = \pi_W(\tau)$. Then we have: 
$$ \Gamma_0 = \{ \widehat{\epsilon}_i - \widehat{\epsilon}_j \ | \ 1 \leq i,j \leq m\}\ \mbox{ and } \ 
\Gamma_s = \{ \widehat{\tau} - \widehat{\epsilon}_i\ | \ 1 \leq i \leq m\}. $$

{\bf In type $S$:} we have $r = m-1$. Let $\pi_S \colon \ZZ^m \to (\ZZ / p\, \ZZ)^{m-1}$ be given by $\pi_S(b_1,\ldots,b_m) = (b_2 - b_1,\ldots,b_m - b_1) \mod p$. Then we have $\pi_S(\epsilon_i) = \widehat{\epsilon}_i$ for $2 \leq i \leq m$. We still denote $\widehat{\epsilon}_1 = \pi_S(\epsilon_1)$, so that actually $\widehat{\epsilon}_1 + \ldots + \widehat{\epsilon}_m = 0$. Let $\widehat{\tau} = \pi_S(\tau)$. Then we have:
$$ \Gamma_0 = \{ \widehat{\epsilon}_i - \widehat{\epsilon}_j \ | \ 1 \leq i,j \leq m\}\ \mbox{ and }\
\Gamma_s = \{ \widehat{\tau} - \widehat{\epsilon}_i - \widehat{\epsilon}_j \ | \ 1 \leq i \neq j \leq m\}. $$

{\bf In type $H$:} we have $r = m / 2$. Let $\pi_H \colon \ZZ^m \to (\ZZ / p\, \ZZ)^{r}$ be given by $\pi_H(b_1,\ldots,b_{2r}) = (b_{r+1} - b_1,\ldots,b_{2r} - b_r) \mod p$. Then we have $\pi_H(\epsilon_i) = \widehat{\epsilon}_i$ for $1 \leq i \leq r$. We still denote $\widehat{\epsilon}_i = \pi_H(\epsilon_i)$ for $i > r$, so that $\widehat{\epsilon}_i +\widehat{\epsilon}_{i+r} = 0$ for $1 \leq i \leq r$ . Let $\widehat{\tau} = \pi_H(\tau)$. Then we have:
$$\Gamma_0 = \{ \widehat{\epsilon}_i + \widehat{\epsilon}_j \ | \ 1 \leq i,j \leq 2r \}\ \mbox{ and }\ 
\Gamma_s = \{ \widehat{\tau} - \widehat{\epsilon}_i  \ | \ 1 \leq i \leq 2r \}. $$

With these explicit descriptions, one can check easily that $\Gamma_0 \cap \Gamma_s = \emptyset$.

\subsubsection{} 
\begin{thm} \label{GeneratingSetsGradedCase}
Let $\FF$ be an algebraically closed field of characteristic $p > 3$. Let $L$ be a simple graded Lie algebra of Cartan type. Then $L$ is generated by 2 elements. \\
\end{thm}

\proof\ Let $X \in \{ W,S,H,K\}$ and $L = X(m,\nn)^{(\infty)}$. In the standard grading, we have $L = L_{-r} \oplus \ldots \oplus L_s$, where $r = 1$ or $2$ and $s > 0$. Furthermore, $L_0$ is isomorphic to $\gl_m$, $\sl_m$, $\sp_{2r}$, $\sp_{2r} \oplus \FF\, I_{2r}$ for  $X = W,S,H,K$ respectively.

Let us first choose an element $y_0 \in L_0$ in an appropriate way. Let $T_X \subseteq L_0$ be the $\ad$-diagonalisable subalgebra given in \ref{torus}. In the natural realisation of $L_0$ as matrices, $T_X$ consists of diagonal matrices. So, we can choose $y_0 \in T_X$ in such a way that the $p$-iterations $y_0,y_0^{[p]},\ldots,y_0^{[p^{n-1}]}$, computed in the restricted Lie algebra $L_0$, span the subspace $T_X$. For example, if $\{ h_j \}_{j \in J}$ is a basis of $T_X$ such that $h_j^{[p]} = h_j$ for all $j$, we can choose $\ds y_0 = \sum_{j \in J} \lambda_j h_j$ such that the coefficients $\{\lambda_j\}_{j \in J}$ are linearly independent over the prime field $\FF_p$.

Now, any Lie subalgebra $L'  \subseteq L$ containing $y_0$ satisfies $[T_X , L'] \subseteq L'$. Indeed, each element $t \in T_X$ can be written as $t = \sum_k \mu_k\, y_0^{[p^k]}$, so that 
$$\ad\, t = \ad \left( \sum_k \mu_k\, y_0^{[p^k]} \right) = \sum_k \mu_k\, (\ad\, y_0)^{p^k}  \in \FF[\ad\, y_0].$$
A priori, the above identity holds as derivations of $L_0$, but Lemma \ref{restrictedaction} ensures that it is also true as derivations of the full Lie algebra $L$. Now, since $L'$ is stable under $\ad\, y_0$, it is also stable under $\ad\, t$.

Let us now construct an element $x \in L$ such that $\FF\, \langle x,y_0 \rangle = L$. We will choose $x$ in the form $x = x_{-1} + x_0 + x_s$, where $x_k \in L_k$ for $k=-1,0$ and $s$. First, by Theorem \ref{GeneratingSetsClassicalCase} or Corollary \ref{corvaria}, there exists $x_0 \in L_0$ such that $\FF\, \langle x_0,y_0 \rangle = L_0$. Now recall that each subspace $L_k$ can be decomposed into weight spaces relative to the action of $T_X$ (see Section \ref{torus}):
$$L_k = \bigoplus_{\alpha \in \Gamma_k} L_k^\alpha.$$
We choose $x_s$ and $x_{-1}$ to be non-zero weight vectors relative to the action of $T_X$. Let $\alpha_s$, $\alpha_{-1}$ be the corresponding weights: by Lemma \ref{torus}, \nr 2, we may assume $\alpha_{-1} \neq \alpha_s$.

Set $L' = \FF \langle x , y_0 \rangle \subseteq L$; we will show that $L' = L$. Let us first decompose $x_0$ according to the weight space decomposition of $L_0$:
$$x_0 = \sum_{\alpha \in \Gamma_0} x_0^\alpha,$$
where each $x_0^\alpha$ is a weight vector associated with the weight $\alpha \in \Gamma_0$. According to Lemma \ref{torus}, we know that none of the weights $\alpha \in \Gamma_0$ is equal to $\alpha_{-1}$ or $\alpha_s$; thus, by our previous choice, the weights $\alpha_{-1}$, $\alpha_s$ and $\alpha \in \Gamma_0$ are pairwise distinct. 

Now we set $\Phi = \Gamma_0 \cup \{ \alpha_{-1} \} \cup \{ \alpha_s \} \subseteq T_X^*$. By Lemma \ref{zariskiopeninh}, there exists $t \in T_X$ such that the values $\left\{ \phi(t) \  | \ \phi \in \Phi \right\} $ are pairwise distinct. Since $y_0 \in L'$, we also have $[T_X,L'] \subseteq L'$. In particular, $(\ad\, t)^k(x) \in L'$ for all $k \geq 0$. Thus, we get:
$$(\ad\, t)^k(x) = \sum_{\phi \in \Phi} \phi(t)^k x^\phi \in L'.$$
As in the proof of Theorem \ref{GeneratingSetsClassicalCase}, this defines an invertible Van der Monde system, so that all weight vectors $x^\phi \in L'$. In particular, $x_{-1}, x_s \in L'$ and $\ds x_0 = \sum_{\alpha \in \Gamma_0} x_0^\alpha \in L'$.

Consequently, we get $L' \supseteq L_0 = \FF \langle x_0 , y_0 \rangle$. Using Proposition \ref{rappels}, \nr 2, we have $\sum_{k \geq 0}(\ad\, L_0)^k (\FF\, x_{-1}) = L_{-1}$ and $\sum_{k \geq 0}(\ad\, L_0)^k (\FF\, x_s) = L_s$, so that also $L' \supseteq L_{-1} + L_s$. Using \nr 3 of the same proposition, we conclude that $L' \supseteq L$, thus $L' = L$. The theorem is proved.

\subsubsection{\bf Remark} Theorem \ref{GeneratingSetsGradedCase} is also true in characteristic $p = 5$ for the Melikian Lie algebras. Let us outline the argument; all required properties of Melikian Lie algebras can be found in \cite[Section 4.3]{Strade}.

Let $M = \cM(n_1,n_2)$ be the Melikian algebra of parameter $(n_1,n_2)$. There exists a grading $M = M_{-3} \oplus \ldots \oplus M_s$, with $s = 3\, (5^{n_1}+5^{n_2})-7$. Also, one has $M_0 \simeq \gl_2$. Let $T \subseteq M_0$ be the natural Cartan subalgebra of $M_0$, corresponding to diagonal matrices in $\gl_2$. Then $T$ acts diagonally on $M$. Let $\Gamma_k$ denote the set of weights of $T$ in the homogeneous part $M_k$. One can see using \cite[relations (4.3.1)]{Strade} that the identity element $I_2 \in \gl_2 \simeq M$ acts by twice the degree, so that $\Gamma_0 \cap \Gamma_{-1} = \Gamma_0 \cap \Gamma_s = \Gamma_{-1} \cap \Gamma_0 = \emptyset$. With these properties, one can proceed as in the proof of Theorem \ref{GeneratingSetsGradedCase} to construct pairs of elements which generate $M$.

\subsection{Generators of simple filtered Lie algebras of Cartan type} $ $\\

In this section we will show that simple filtered algebras of Cartan type are generated by 2 elements. The method is essentially a reduction to the graded case. However, due to the fact that the elements constructed in Theorem \ref{GeneratingSetsGradedCase} are not homogeneous, usual filtered algebra arguments do not apply in any straightforward manner.

We first need some general results about linear algebra in characteristic $p$. The underlying idea is to replace polynomials in endomorphisms by $p$-polynomials. This will enable us to construct tori acting on filtered Lie algebras in a way which is compatible with the filtration.

\subsubsection{\bf Additive polynomials and $p$-polynomials}
A {\it $p$-polynomial} is a polynomial of the form $f(t) = a_0\, t + a_1\, t^p + \ldots + a_n\, t^{p^n} \in \FF[t]$. From \cite[Lemmas 20.3.A and 20.3.B]{Humphreys}, we have: \\

\begin{lemma} Let $f \in \FF[t]$ be any polynomial.
\begin{enumerate}
\item The polynomial $f$ is additive if and only if $f$ is a $p$-polynomial.
\item If $f$ is a $p$-polynomial, then the roots of $f$ form a subgroup of $\FF$.
\item Conversely, if $f$ is square-free and its roots form a subgroup of $\FF$, then $f$ is a $p$-polynomial.
\end{enumerate}
\end{lemma}

\subsubsection{} Let $V$ be a finite-dimensional vector space over $\FF$ and $u \in \End(V)$ a linear map. Denote by $\egv(u)$ the set of eigenvalues of $u$, and $\Lambda(u)$ the additive subgroup of $\FF$ generated by $\egv(u)$. Since the additive subgroups of $\FF$ are the sub-vector spaces over the prime field $\FF_p$, we can define the {\it $p$-order of $u$} to be $\ord(u) = \dim_{\FF_p}\, \Lambda(u)$; alternatively, it is the integer $r \geq 0$ such that $| \Lambda(u) | = p^r$.

Also, there exists a $p$-polynomial $f(t)$ of lowest degree such that $f(u) = 0$. It is unique up to a scaling factor; we call it the {\it minimal $p$-polynomial of $u$}. \\

\begin{lemma} \label{propertiesporder}
Let $u \in \End(V)$ and $\displaystyle f_u(t) = \prod_{\lambda \in \Lambda(u)} (t-\lambda)$.
\begin{enumerate}
\item For all $k \geq 0$, we have $\ord(u^{p^k}) = \ord(u)$.
\item Assume that $u$ is semi-simple. Then $f_u(t)$ is the minimal $p$-polynomial of $u$. In particular, let $r = \ord(u)$: then $\deg(f_u) = p^r$. Also, if $T = \sum_{j \geq 0} \FF\, u^{p^j}$, then $\dim(T) = r$.
\item There exists a minimal integer $k \geq 0$ such that $u^{p^k}$ is semi-simple. The minimal $p$-polynomial  of $u$ is $(f_u(t))^{p^k}$. \\
\end{enumerate}
\end{lemma}

\proof\ For \nr 1, we show the result for $k=1$; the general case follows by induction. We consider the Jordan decomposition $u = s+n$, where $s$ is semi-simple, $n$ nilpotent and $[s,n]=0$. Note that $u^p = s^p + n^p$ is the Jordan decomposition of $u^p$.

Since $u$ and $s$ have the same eigenvalues, we have $\Lambda(u) = \Lambda(s)$. Also, $s$ being semi-simple, it is easy to see that the eigenvalues of $s^p$ are the $p$-th powers of eigenvalues of $s$: thus, the Frobenius automorphism of $\FF$ induces an isomorphism $\Lambda(s) \stackrel{\sim}{\to} \Lambda(s^p)$. Hence, we have:
$$\Lambda(u^p) = \Lambda(s^p) \simeq \Lambda(s) = \Lambda(u);$$
it immediately follows that $\ord(u^p) = \ord(u)$.

Let us prove \nr 2. By semi-simplicity of $u$, we know that a polynomial $g(t)$ annihilates $u$ if and only if the eigenvalues of $u$ are roots of $g(t)$. Thus, a $p$-polynomial annihilates $u$ if and only if its set of roots contains $\Lambda(u)$. The fact that $f_u(t)$ is the minimal $p$-polynomial of $u$ follows easily. Also, by construction we have $\deg(f_u) = |\Lambda(u)| = p^r$, where $r = \ord(u)$.

Last, let $T = \sum_{j \geq 0} \FF\, u^{p^j}$. Since $u$ is annihilated by a $p$-polynomial of degree $p^r$, it is clear that $\dim(T) \leq p^r$. Conversely, if the elements $u,u^p,\ldots,u^{p^{r-1}}$ were not linearly independent, then $u$ would be annihilated by a $p$-polynomial of degree $\leq p^{r-1}$, a contradiction: hence $\dim(T) \geq r$. 

Let us prove \nr 3. Again we consider the Jordan decomposition $u = s + n$. Since $u^{p^k} = s^{p^k} + n^{p^k}$ is the Jordan decomposition of $u^{p^k}$, this endomorphism is semi-simple if and only if $n^{p^k} = 0$. By nilpotence of $n$, there exists a minimal value of $k$ satisfying this property. Since the eigenvalues of $u^{p^k}$ are the $p^k$-th powers of the eigenvalues of $u$, the minimal $p$-polynomial of $u^{p^k}$ is (by \nr 2):
$$g(t) = f_{u^{p^k}}(t) = \prod_{\lambda \in \Lambda(u)}(t - \lambda^{p^k}).$$
By definition, $g (u^{p^k}) = 0$. Thus, the polynomial $g(t^{p^k})$ annihilates $u$; but we have:
$$g(t^{p^k}) = \prod_{\lambda \in \Lambda(u)}(t^{p^k} - \lambda^{p^k}) = \prod_{\lambda \in \Lambda(u)}(t - \lambda) ^{p^k} = f_u(t)^{p^k}.$$
To prove that $f_u(t)^{p^k}$ is minimal, one can proceed as follows. First, we show that the minimal $p$-polynomial has the form $F(t) = f_u(t)^{p^s}$ for some integer $s$: indeed, the roots of $F$ are exactly the elements of $\Lambda(u)$, so that $F$ has the form:
$$F(t) = \prod_{\lambda \in \Lambda(u)} (t-\lambda)^{m_\lambda},$$
and because $F$ is a $p$-polynomial, the multiplicities $m_\lambda$ have to be some power of $p$ and independent of $\lambda$. Then we argue as before that the polynomial: 
$$G(t) = \prod_{\lambda \in \Lambda(u)} (t-\lambda^{p^s})$$
annihilates $u^{p^s}$. Since $G(t)$ has only simple roots, this implies that $u^{p^s}$ is semi-simple. By minimality of $k$, we get $s \geq k$, proving minimality of $f_u(t)^{p^k}$. The lemma is proved.

\subsubsection{\bf Filtered algebras}
We recall shortly some general results about filtered Lie algebras. Let $L$ be a filtered Lie algebra. We assume the filtration to be descending, i.e. $L = L_{(-s')} \supseteq \ldots \supseteq L_{(s)} \supseteq 0$, with $s,s' \geq 0$. Let $\cG = Gr(L)$ be the associated graded Lie algebra. Homogeneous components in $\cG$ are denoted by $\cG_k$, for $k \in \{ -s',\ldots,s\}$. Also, denote by $gr_k$ the natural linear map $gr_k : L_{(k)} \twoheadrightarrow \cG_k = L_{(k)} / L_{(k+1)}$.

For any $x \in L \smallsetminus 0$, we denote by $\nu(x)$ the largest integer $k$ such that $x \in L_{(k)}$. We set $gr(x) = gr_k(x) = x + L_{(k+1)} \in Gr(L)$. If $x = 0$ we denote $gr(x) = 0$. By construction, we have $gr(x) \neq 0$ if $x \neq 0$, and $\deg\, gr(x) = \nu(x)$.

Denote by $\{ .,. \}$ the Lie bracket in $Gr(L)$. It is uniquely defined by its values on homogeneous elements. Let $x,y \in L$ be non-zero elements and $a = \nu(x)$, $b = \nu(y)$. Then $\{ gr(x) , gr(y) \} = gr_{a+b}\, \big( [x,y] \big) = [x,y] + L_{(a+b+1)}$. By construction, if $\{ gr(x) , gr(y) \} \neq 0$, then $\{ gr(x) , gr(y) \} = gr\, \big( [x,y] \big)$.

\subsubsection{}
Let $x \in L$. We can apply results of \ref{propertiesporder} to the endomorphism $\ad\, x \in \End(L)$. The set of eigenvalues and the generated subgroup will be denoted by $\egv(x)$ and $\Lambda(x)$ instead of $\egv(\ad\, x)$ and $\Lambda(\ad\, x)$. Also, define the {\it $p$-order of $x$} to be $\ord(x) = \ord(\ad\, x)$. \\

\begin{lemma} \label{ordh=ordgrh}
Let $L$ be a filtered Lie algebra and $h \in L_{(0)}$. Then $\ord(h) = \ord(gr\, h)$. \\
\end{lemma}

\proof\ If $h \in L_{(1)}$, then both $h$ and $gr\, h$ are $\ad$-nilpotent, so that their $p$-orders are 0. Now we assume that $h \in L_{(0)} \smallsetminus L_{(1)}$, i.e. $\nu(h) = 0$. First we establish that $\egv(h) \subseteq \egv(gr\, h)$, so that $\ord(h) \leq \ord (gr\, h)$. Let $\lambda \in \FF$ be an eigenvalue of $\ad\, h$: there exists $x \in L$ such that $[h,x] = \lambda\, x$. Let $b = \nu(x)$. Since $\nu(h) = 0$, we obtain:
$$\{ gr(h) , gr(x) \} = \lambda\, x + L_{(b + 1)} = \lambda\, gr(x),$$
so that $\lambda \in \egv(gr\, h)$.

Let us show now that $\egv(gr\, h) \subseteq \Lambda(h)$. This would imply $\Lambda(gr\, h) \subseteq \Lambda(h)$, and therefore $\ord(gr\, h) \leq \ord(h)$ as desired. Let $\lambda \in \egv(gr\, h)$. There exists $\xi \in Gr(L)$ such that $\{ gr(h) , \xi \} = \lambda\, \xi$. Since the Lie bracket $\{ .,. \}$ is homogeneous, we may assume that $\xi = gr(x)$ for some $x \in L$. Let $b = \nu(x)$: the assumption means that there exists $y \in L_{(b+1)}$ such that $[h,x] = \lambda\, x + y$.

Let $\varphi$ be the minimal $p$-polynomial annihilating $u = \ad(h) \in \End(L)$. In particular, the roots of $\varphi$ are exactly the elements of $\Lambda(h)$. Let $\psi$ be the (ordinary) polynomial such that $\varphi(t) = t\, \psi(t)$. We have $(u - \lambda\, \id)\, x = y$. By additivity of $\varphi$, we get:
\begin{equation} \label{lhs}
\varphi(u - \lambda\, \id)\, x = \varphi(u)\, x - \varphi(\lambda)\, x = -\varphi(\lambda)\, x.
\end{equation}
Also, using $\varphi(t) = \psi(t)\, t$, we get:
\begin{equation} \label{rhs}
\varphi(u - \lambda\, \id)\, x = \psi(u-\lambda\, \id)\, (u - \lambda\, \id)\, x = \psi(u-\lambda\, \id)\, y.
\end{equation}
Since $u = \ad(h)$ with $h \in L_{(0)}$, it is easy to see that $\psi(u-\lambda\, \id) \in \End(L)$ preserves the filtration. By (\ref{rhs}), we get $\nu(\psi(u - \lambda\, \id)\, y) \geq \nu(y) > b = \nu(x)$. This is compatible with (\ref{lhs}) only if $\varphi(\lambda) = 0$, i.e. $\lambda \in \Lambda(h)$. The lemma is proved. \\

{\nbf Remark.} The result doesn't hold if $\nu(h) < 0$. For instance, let $L$ be the first Witt algebra $L = W(1,1) = \FF\, e_{-1} \oplus \ldots \oplus \FF\, e_{p-2}$, with bracket $[e_i,e_j] = (j-i)e_{i+j}$ when defined. We equip $L$ with the natural filtration $L_{(k)} = \FF\, e_k \oplus \ldots \FF\, e_{p-2}$. Let $h = e_{-1}+e_0$. Then $h$ is semi-simple, with eigevalues $0,1,\ldots,p-1 \in \FF$ while $gr(h) = e_{-1} + L_{(0)}$ is nilpotent. Thus, $ \ord(h) = 1 \neq 0 = \ord\, gr(h)$.

\subsubsection{} \label{GeneratingSetsFilteredCase}
\begin{thm}
Let $\FF$ be an algebraically closed field of characteristic $p > 3$ and $L$ be a simple filtered Lie algebra of Cartan type over $\FF$. Then $L$ is generated by 2 elements. \\
\end{thm}

\proof\ We denote $\cG = Gr(L)$ the associated graded algebra of $L$, so that the last derived subalgebra $\cG^{(\infty)}$ is simple graded of Cartan type. Let $\eta \in \cG^{(\infty)}_0$ be regular semi-simple as in the proof of Theorem \ref{GeneratingSetsGradedCase}. There exists $y \in L$ such that $gr(y) = \eta$. For sufficiently large $k$, the linear map $(\ad\, y)^{p^k} \in \End(L)$ is semi-simple. We let:
$$T_L = \Span_\FF \{ (\ad\, y)^{p^k},(\ad\, y)^{p^{k+1}},\ldots\} \subseteq \End(L)$$
and
$$T_{\cG} = \Span_\FF \{ (\ad\, gr\, y)^{p^k},(\ad\, gr\, y)^{p^{k+1}}, \ldots\} \subseteq \End(\cG).$$
By Lemmas \ref{ordh=ordgrh} and \ref{propertiesporder}, $\dim(T_L) = \dim(T_\cG)$. Let $\gamma:T_L \to T_\cG$ be the linear isomorphism defined by $(\ad\, y)^{p^j} \mapsto (\ad\, gr\, y)^{p^j}$. We will often identify $T_L$ and $T_\cG$ via $\gamma$ in the sequel.

By semi-simplicity of $(\ad\, y)^{p^k}$, it is easy to see that $T_L$ acts semi-simply on $L$. Since $y \in L_{(0)}$, it follows that $T_L$ respects the filtration of $L$. Thus, for all $j$, there exists a $T$-stable subspace $V_j \subseteq L$ such that $L_{(j)} = V_j \oplus L_{(j+1)}$. Recall the natural projection $gr_j:L_{(j)} \twoheadrightarrow \cG_j = L_{(j)} / L_{(j+1)}$. Then $gr_j$ restricts to an isomorphism of vector spaces $gr_j: V_j \to \cG_j$. This isomorphism is actually an isomorphism of $T_L$-modules. This means that the following diagram is commutative for all $h \in T_L$:
$$\begin{CD}
V_j   @>{gr_j}>> \cG_j \\
@V{h}VV \hspace{-1cm} \circlearrowleft \hspace{.5cm} @VV{ \gamma(h)}V\\ 
V_j                  @>>{gr_j}> \cG_j
\end{CD}$$
Commutativity of the diagram is then easy to check when $h = (\ad\, y)^{p^i}$ for some $i$.

Let $\Gamma(T_L,V_j)$ (resp. $\Gamma(T_\cG,\cG_j)$) be the set of roots of $T_L$ acting on $V_j$ (resp. of $T_\cG$ acting on $\cG_j$). If we identify $T_L$ and $T_\cG$ via the isomorphism $\gamma$, then by the previous discussion, we have $\Gamma(T_L,V_j) = \Gamma(T_\cG,\cG_j)$.

Now we are ready to construct pairs of generators of the Lie algebra $L$. Let $t$ be the maximal degree of homogeneous components of the simple Lie algebra $\cG^{(\infty)}$. We keep $y$ as above; we will construct an element $x = x_{-1}+x_0+x_t \in L$ similarly to what was done in the proof of Theorem \ref{GeneratingSetsGradedCase}. 

We first choose three elements $\xi_{i} \in \cG^{(\infty)}_{i}$, for $i \in \{ -1,0,t\}$, in such a way that $\FF\langle \eta,\xi_0 \rangle = \cG^{(\infty)}_0$ and $\xi_{-1},\xi_t$ are weight vectors for $T_\cG$ associated with different weights. For $i \in \{ -1,0,t \}$, let $x_i \in V_i$ be such that $\xi_i = gr(x_i)$; we can arrange $x_{-1}$ and $x_t$ to be weight vectors. By the previous discussion and Lemma \ref{torus}, the weights involved in the decompositions of $x_{-1}$, $x_0$ and $x_t$ are the same as for $\xi_{-1}$, $\xi_0$ and $\xi_t$ and therefore pairwise distinct.

Now set $x = x_{-1} + x_0 + x_t$ and $M = \FF \langle y , x \rangle \subseteq L$. Since $M$ is stable under $(\ad\, y)$, it is also $\ad T_L$-stable; arguing as in the proof of Theorem \ref{GeneratingSetsGradedCase}, we obtain that $x_{-1}$, $x_0$, $x_t \in M$. Thus, we have $\xi_i = gr(x_i) \in Gr(M)$ for $i \in \{ -1,0,t\}$ and $\eta = gr(y) \in Gr(M)$. By the construction of Theorem \ref{GeneratingSetsGradedCase}, we have $\cG^{(\infty)} = \FF \langle \xi_{-1}+\xi_0+\xi_t, \eta \rangle \subseteq Gr(M)$. Since $Gr(M) \subseteq Gr(L) = \cG$, it follows that $M \subseteq L$ satisfies: 
$$(Gr\, M)^{(\infty)} = (Gr\, L)^{(\infty)} \neq 0.$$
By \cite[Lemma 4.2.5]{Strade} this implies that $M$ contains a non-zero ideal of $L$. By simplicity of $L$, we must have $M = L$ as desired.

\subsection{One-and-a-half generation in $W(m,\nn)$} $ $ \\

In this section we examine the question whether the Lie algebra $W(m,\nn)$ is generated by one and a half elements. As it turns out, the answer is positive if and only if $m = 1$. In view of \ref{plongement}, it is natural to first carefully examine the case $\nn = 1$.

\subsubsection{\bf Notations} \label{debutnotationW}
Let $m \geq 1$ be a positive integer. As in \ref{lanotationBm}, we work with the truncated polynomial ring $B_m$ instead of the divided power algebra $\cO(m,\un)$. For a multi-index $\alpha \in \NN^m$, we will write $x^\alpha = x_1^{\alpha_1} \cdots  x_m^{\alpha_m}$, so that in particular there is a unique non-zero monomial of highest degree, say $x^\tau$ for $\tau = (p-1,\ldots,p-1)$. We have $\deg(x^\tau) = m(p-1)$.

The Lie algebra $W(m,\un) = \Der\, \cO(m,\un)$ is simple, restricted and graded. We denote $W(m,\un) = W_{-1} \oplus \ldots \oplus W_s$, where:
$$ W_{k} = \bigoplus_{j=1}^m \sum_{|\alpha| = k+1} \FF\, x^\alpha D_j.$$
We have $s = m(p-1)-1$ and $\dim\, W(m,\un)_s = \dim W(m,\un)_{-1} = m$.

In the sequel, we will usually write $W$ instead of $W(m,\un)$ and $B$ instead of $B_m$. The Lie algebra $W$ is also a $B$-module, so that products $b \, w$ for $b \in B$ and $w \in W$ make sense. Note that any homogeneous element $\delta \in W_s$ of maximal degree  can be uniquely written as $\delta = x^\tau\, \d$, for some non-zero derivation $\d \in W_{-1}$.

\subsubsection{} \label{introy}
Consider any element $y \in W$. First we gather some information on $y$ and its action on $W$. For all $j \geq 0$, we denote by $y^{[p^j]}$ the unique element of $W$ satisfying $(\ad\, y)^{p^j} = \ad(y^{[p^j]})$. For any $p$-polynomial $f(t) = \sum \alpha_j \, t^{p^j}$, we can define an element $f(y) \in W$ by setting $f(y) = \sum \alpha_j\, y^{[p^j]}$. It is easy to see using the Leibniz rule that $\ad\, f(y) = f(\ad\, y) \in \Der\, (W)$.

For all $j \geq 0$, there exists $\d_j \in W_{-1}$ such that:
$$y^{[p^j]} \equiv \d_j \mod W_{\geq 0}.$$
We will denote by $k \geq 0$ the least integer such that the elements $\d_0,\ldots,\d_k \in W_{-1}$ are linearly dependent. We may have $k = 0$; in any case $k \leq m$. Let $\alpha_0\, \d_0 + \ldots + \alpha_k\, \d_k =0$ be a nontrivial linear dependence relation; in particular $\alpha_k \neq 0$. Let $\ds f(t) = \sum_{j=0}^k \alpha_j \, t^{p^j}$: by construction, $f(y) \in W_{\geq 0}$. Let $h \in W_0$ be the degree 0 component of $f(y)$, so that $f(y) \equiv h \mod W_{>0}$.

Since $f(t)$ has no constant term, there exists $g(t) \in \FF[t]$ such that $f(t) = t\, g(t)$. Note that $g(t)$ is not a $p$-polynomial, so that $g(\ad\, y) \in \End\, (W)$ is still defined as an endomorphism of $W$, but is not a derivation anymore in general.

\subsubsection{} \label{generalitesady}
Recall that for any integer $a \geq 0$, one has a $p$-ary expansion $a = a_0 + p\, a_1 + \ldots + p^l a_l$, with $0 \leq a_i < p$ for all $i$. We will write $a = [a_l,\ldots,a_0]$. The decomposition is unique up to some zero coefficients to the left. We also define the $p$-ary length of $a$ to be the sum of its $p$-ary digits, i.e. $|a|_p = a_l + \ldots + a_0$. Properties of $p$-ary expansions are often intuitive if one thinks of the usual decimal representation of integers. \\

\begin{lemma} 
Let $y\in W,\ k\in \NN,\ g(t)\in \FF[t]$ be as in \ref{introy}. Let $x \in W_s$, $x \neq 0$.
\begin{enumerate}
\item For any integer $a \geq 0$, one has $(\ad\, y)^a(x) \in W_{\geq s - |a|_p}$. 
\item If $a < p^k$, then $(\ad\, y)^a(x) \not \in W_{> s - |a|_p}$.
\item Let $\delta = g(\ad\, y)(x) \in W$. Then $\delta \in W_{s-k(p-1)} \smallsetminus W_{>s-k(p-1)}$. Furthermore, $x \in B\, \delta$. \\
\end{enumerate}
\end{lemma}

\proof\ Let $a = [a_{l},\ldots,a_0]$ be the $p$-ary expansion of $a$. We denote $q = |a|_p$ the $p$-ary length of $a$. Then we have:
$$(\ad\, y)^a = (\ad\, y)^{a_0 + p a_1 + \ldots + p^l a_l} = (\ad\, y)^{a_0} \ldots (\ad\, y^{[p^l]})^{a_l}.$$
Since $[W,W_{\geq d}] \subseteq W_{\geq d-1}$ for all $d$, we immediately get:
$$(\ad\, y)^a(W_s) \subseteq W_{\geq s-a_0-\ldots-a_l} = W_{\geq s - q}.$$
Recall that $x$ can be written as $x = x^\tau \d$ with $\d \in W_{-1} \smallsetminus \{0\}$.  Furthermore, recall that for all $j \geq 0$, we have $y^{[p^j]} = \d_j + r_j$ with $\d_j \in W_{-1}$ and $r_j \in W_{\geq 0}$. Thus,
\begin{equation*}
(\ad\, y)^a = (\ad\, \d_0 + \ad\, r_0)^{a_0} \cdots (\ad\, \d_l + \ad\, r_l)^{a_l}  = (\ad\, \d_0)^{a_0} \cdots (\ad\, \d_l)^{a_l} + \ldots,
\end{equation*}
where ``$\ldots$'' is a sum of products of the form $(\ad\, u_1) \cdots (\ad\, u_q)$, with $u_i \in \{ \d_i, r_i \}$ for all $i$, but not all $u_i = \d_i$. Since one has $[W_0,W_{\geq d}] \subseteq W_{\geq d}$ for all $d$, we can see easily that for such a product:
$$(\ad\, u_1) \cdots (\ad\, u_q)(W_s) \subseteq W_{> s-q}.$$
We can deduce:
\begin{equation} \label{intermediaire}
(\ad\, y)^a (x) \equiv (\ad\, \d_0)^{a_0} \cdots (\ad\, \d_l)^{a_l}(x) \mod W_{>s - q}.
\end{equation}
Now, recall the identity, for any $D,E \in W$ and $b \in B$:
$$(\ad\, D)(bE) = [D,bE] = D(b) E + b [D,E].$$
Using this identity and the fact that elements of $W_{-1}$ commute to each other, equation (\ref{intermediaire}) becomes:
\begin{equation} \label{intermediaire2}
(\ad\, y)^a (x) \equiv \d_0^{a_0} \cdots \d_l^{a_l}(x^\tau) \cdot \d \mod W_{>s - q}.
\end{equation}
Now, since $a < p^k$, we can see that the $p$-ary expansion $a = [a_l,\ldots,a_0]$ has at most $k$ digits, i.e. $l \leq k-1$. Thus, allowing the $p$-ary expansion to start with some $0$'s, we can also write $a = [a_{k-1},\ldots,a_0]$. By definition of $k$, the elements $\d_0,\ldots,\d_{k-1}$ are linearly independent, so that all coefficients $\d_0^{a_0} \cdots \d_{k-1}^{a_{k-1}}(x^\tau) \neq 0$ by Lemma \ref{usefulofwm1}. Since the action of $B$ on $W_{-1}$ is faithful, we also get $\d_0^{a_0} \cdots \d_{k-1}^{a_{k-1}}(x^\tau) \cdot \d \neq 0$, so that $(\ad\, y)^a (x) \not \in W_{>s - q}$. This proves claim 2.

Now for claim 3 we just have to make the following observation. The element $\delta = g(\ad\, y)(x)$ is a sum of terms of the form $(\ad\, y)^{p^j-1}(x)$ with $j \leq k$, so that according to \nr 2, we have $(\ad\, y)^{p^j-1}(x) \in W_{s-|p^j-1|_p} \smallsetminus W_{> s-|p^j-1|_p}$. One can easily compute that $|p^j-1|_p = j(p-1)$. Since the term of degree $p^k-1$ in $g(t)$ is non-zero, it follows that $g(\ad\, y)(x) \in W_{s-k(p-1)} \smallsetminus W_{> s-k(p-1)}$. Using relation (\ref{intermediaire2}), we see more precisely that there exists a homogeneous element $\beta \in B$ such that $g(\ad\, y)(x) \equiv \beta \d \mod W_{> s-k(p-1)}$. By Proposition \ref{usefulofomn}, there exists $b \in B$ such that $b\, \beta = x^\tau$; thus, $b\cdot \delta \equiv b\, \beta\cdot \d \mod W_{>s} = x^\tau\,\d + \{0\}$, so that $x = b\cdot \delta \in B\, \delta$ as required. The lemma is proved.

\subsubsection{\bf The case $k = m$} We keep the notations $W = W(m,\un)$ and $B = \cO(m,\un)$ of \ref{debutnotationW}. \\

\begin{prop} \label{casek=m}
Let $x \in W_s$ and $y \in W$ be non-zero elements. Let $L = \FF \langle x,y \rangle \subseteq W$ be the subalgebra generated by $x$ and $y$. Define $k \in \NN$ and $g(t) \in \FF[t]$ as in \ref{introy} and assume that $k = m$. Let $\delta = g(\ad\, y)(x) \in W$. Then $[y,\delta] = 0$ and $L \subseteq B\, \delta + \FF\, y$. \\
\end{prop}

\proof\ Let us recall from \ref{introy} the construction of $g(t)$. The assumption $k = m$ means that the elements $\d_0,\ldots,\d_{m-1}$ are linearly independent. Then we take some non-trivial relation $\alpha_0\, \d_0 + \ldots + \alpha_m\, \d_m = 0$; observe that the coefficients $\alpha_0,\ldots,\alpha_m$ are uniquely determined up to some non-zero scalar factor. Set $f(t) = \alpha_0\, t + \ldots + \alpha_m\, t^{p^m}$ and let $g(t)$ be such that $f(t) = t g(t)$. In particular, we have $f(\ad\, y) = (\ad\, y)\circ g(\ad\, y)$.

Using \cite[Corollary 1]{Premet_Wn}, we can see that there exists a non-zero $p$-polynomial of degree at most $m$, say $F(t) = \beta_0 t + \ldots + \beta_m t^{p^m}$ which annihilates $\ad\, y$. We stress the fact that the integer $m$ is the one defining $W = W(m,\un)$. We have $F(y) = \beta_0\, y + \ldots + \beta_m\, y^{[p^m]} = 0$. Taking the homogeneous part of degree $-1$, we get $\beta_0\, \d_0 + \ldots + \beta_m\, \d_m = 0$. By the previous observation, we see that all coefficients $\beta_j$ are equal to $\alpha_j$ up to some non-zero scalar factor, which we may assume to be 1. Thus, we obtain $f(t) = F(t)$, so that $f(y) = 0$ and hence $f(\ad\, y) = \ad\, f(y) = 0$. It follows:
\begin{equation*}
[ y,\delta ] = (\ad\, y)(\delta) = (\ad\, y) \circ g(\ad\, y)(x) = f(\ad\, y)(x) = 0.
\end{equation*}
Using this relation it is easy to show that the subspace $B\, \delta + \FF\, y$ is actually a Lie subalgebra of $W$. Using Lemma \ref{generalitesady}, \nr 3, we see that $x \in B\, \delta$, so that the Lie algebra $B\, \delta + \FF\, y$ contains both $x$ and $y$. Consequently, we have $L(x,y) \subseteq B\, \delta + \FF\, y$. The proposition is proved.

\subsubsection{\bf The case $k < m$} We keep the notations $W$ and $B$ as above. \\

\begin{prop} \label{casek<m}
Let $x \in W_s$ and $y \in W$ be non-zero elements. Let $L = \FF \langle x,y \rangle \subseteq W$ be the subalgebra generated by $x$ and $y$. Define $k \in \NN$ as in \ref{introy} and assume that $k < m$. Let $H = \sum_{j \geq 0} \FF\, (\ad\, y)^j(x)$. Then $L = H + \FF\, y$. \\
\end{prop}

\proof\ Let us recall some facts from \ref{introy}. There exists a $p$-polynomial $f(t)$ of degree $p^k$ such  that $f(y) \in W_{\geq 0}$. As a consequence, $\ad\, f(y)(W_{\geq d}) \subseteq W_{\geq d}$ for all $d$; in particular, $W_s$ is stable under $\ad\, f(y)$. Hence, the subspace $\ds H_s = \sum_{j \geq 0}(\ad\, f(y))^j(x) \subseteq H \cap W_s$.

Let us show first that $H \subseteq W_{\geq s-k(p-1)}$. For all integers $N \geq 0$, there exists a decomposition $\ds t^N = \sum_{j \geq 0} \varphi_j(t) f(t)^j$, with all coefficients $\varphi_j(t) \in \FF[t]$ of degree less than $\deg(f) = p^k$. Using the identity $f(\ad\, y) = \ad\, f(y)$, we obtain: 
$$(\ad\, y)^N(x) = \sum_{j \geq 0} \varphi_j(\ad\, y) \underbrace{(\ad\, f(y))^j(x)}_{\in H_s \subseteq W_s}.$$
Now all $\varphi_j(\ad\, y)$ are linear combinations of the $(\ad\, y)^i$ with $i < p^k$. Any such integer $i$ has a $p$-ary expansion with $k$ digits, so that the $p$-ary length $|i|_p \leq k(p-1)$. In light of Lemma \ref{generalitesady}, \nr 1, we get:
$$(\ad\, y)^N(x) \in \sum_{i=0}^{p^k-1} (\ad\, y)^i (W_s) \subseteq \sum_{i=0}^{p^k-1} W_{\geq s-|i|_p} \subseteq W_{\geq s - k(p-1)}.$$

Now we check that $[x,H] = 0$. Recall that the maximal degree of a homogeneous derivation is $s = m(p-1) - 1$. Since $k < m$ and $p > 3$, it implies $s - k(p-1) \geq 1$, so that $[x, H] \subseteq [x,W_{\geq s-k(p-1)}] \subseteq [W_s , W_{\geq 1}] = 0$.

Note that by construction $x,y \in H + \FF\, y$. Using the relation $[x,H] = 0$, we can see that $H + \FF\, y$ is stable under $\ad\, x$ and $\ad\, y$. The subalgebra $L(x,y)$ is also the smallest subspace of $W$ containing $x,y$ and stable under $\ad\, x$ and $\ad\, y$, hence $L(x,y) \subseteq H + \FF\, y$. The reverse inclusion is clear: the proposition is proved.

\subsubsection{} \label{cork}

\begin{cor}
Let $W = W(m,\un)$. Let $x \in W_s$ and $y \in W$. Let $L \subseteq W$ be the subalgebra generated by $x$ and $y$. Then $\dim\, [L,L] \leq p^m$. \\
\end{cor}

\proof\ Let $k$ be defined as in \ref{introy}. Denote by $B = \cO(m,\un)$. If $k = m$, then $L \subseteq B\, \delta + \FF\, y$ for some element $\delta \in W$ such that $[\delta,y] = 0$ (Proposition \ref{casek=m}). In this case, we have $[L,L] \subseteq B\, \delta$; thus, $\dim\, [L,L] \leq \dim(B) = p^m$.

If $k < m$, let $H$ be the sub-$\FF\, y$-module generated by $x$. Recall that $L = H + \FF\, y$, with $[x,H] = 0$, $[y,H] \subseteq H$. Since $L$ is generated by $x$ and $y$, it follows readily that $[L,L] = [L,H + \FF\, y] \subseteq H$. Since $\ad\, y$ is annihilated by a polynomial of degree $p^m$ \cite[Corollary 1]{Premet_Wn}, we deduce that $\dim\, H \leq p^m$ as desired.

\subsubsection{\bf Case of the Zassenhaus algebra} \label{zassenhaus}
Now we can turn to the problem of one and a half generation in general graded Lie algebras of type $W$. We will first study the case of the Zassenhaus Lie algebras $W(1,n)$ for arbitrary $n \geq 1$. The method is heavily computational.

The Lie algebra $W(1,n)$ has a decomposition into graded components:
$$W(1,n) = W_{-1} \oplus W_0 \oplus \ldots \oplus W_s,$$
where $s = p^n-2$, each $W_k = \FF\, e_k$ has dimension 1 and the following relation holds:
\begin{equation} \label{relationszassenhaus}
[e_i,e_j] = \left\{ \begin{array}{ll}
\left( \binom{i+j+1}{j} - \binom{i+j+1}{i} \right) e_{i+j} & \mbox{if } -1 \leq i+j \leq s, \\
0 & \mbox{otherwise}.
\end{array} \right.
\end{equation}
In particular, we have $[e_{-1},e_k] = e_{k-1}$ for all index $k$ and:
$$[e_k,e_s] = \left\{ \begin{array}{lcl}
e_{s-1} & \mbox{if} & k=-1\ ;\\
-2\, e_s & \mbox{if} & k=0\ ;\\
0 & \mbox{if} & k > 0.
\end{array} \right.$$

\begin{prop}
Let $x \in W(1,n)$ be a non-zero element. There exists $y \in W(1,n)$ such that $x$ and $y$ generate $W(1,n)$. \\
\end{prop}

\proof\ First we deal with the easy case $x \in W_s$. We can then choose $y = e_{-1}$, so that $\FF\langle x,y \rangle = \FF\langle e_s,e_{-1} \rangle = W(1,n)$. From now on we will assume that $x \not \in W_s$.

We can decompose $x$ in a sum: $ x = \lambda_{-1}\, e_{-1} + \ldots + \lambda_s\, e_s$, for some scalars $\lambda_i \in \FF$; by assumption there exists $j < s$ such that $\lambda_j \neq 0$. We will represent $x$ as a column:
$$x \equiv \left[ \begin{array}{c}
\lambda_{-1}  \\ 
\lambda_0 \\
\vdots \\
\lambda_s 
\end{array} \right].$$
Let $y_\alpha = e_{-1} + \alpha\, e_s$, where $\alpha \in \FF$ is a scalar to be chosen later. We consider $\FF\langle x,y_\alpha \rangle$, the subalgebra generated by $x$ and $y_\alpha$. Let us compute the bracket $[y_\alpha,x]$. We easily obtain: 
$$(\ad\, y)(x) \equiv  \left[ \begin{array}{c}
\lambda_0  \\
\lambda_1 \\
\vdots  \\
\lambda_s - \alpha\, \lambda_{-1} \\
2 \alpha\, \lambda_0
 \end{array} \right].$$
By induction, we find for $k \leq s$:
$$(\ad\, y_\alpha)^k(x) \equiv  \left[ \begin{array}{c}
\lambda_{k-1} \\
\lambda_k \\
\vdots  \\
\lambda_s - \alpha\, \lambda_{-1}  \\
\alpha\, \lambda_0 \\
\vdots \\
\alpha\, \lambda_{k-2} \\
2 \alpha\, \lambda_{k-1}
\end{array} \right].$$

Let us denote $x^{(k)} = (\ad\, y_\alpha)^k(x) \in L(x,y_\alpha)$. We will show that, for a suitable choice of $\alpha$, the family $\{y,x,x^{(1)},\ldots,x^{(s)}\}$ has rank $s+2 = \dim\, W(1,n)$. From this property, it would follow that $\dim\, \FF\langle x,y_\alpha \rangle \geq \dim\, W(1,n)$, forcing $\FF\langle x,y_\alpha \rangle =  W(1,n)$. We need to show that the following matrix $M_\alpha$ has non-zero determinant for some values of $\alpha$:
$$ M_\alpha = \left[ \begin{array}{ccccc}
1 & \lambda_{-1} & \lambda_0 & \ldots & \lambda_{s-1} \\
0 & \lambda_0 & \lambda_1 & \ldots & \lambda_s - \alpha\, \lambda_{-1} \\
\vdots & \vdots & \vdots & \ldots & \alpha\, \lambda_0 \\
\vdots & \vdots & \vdots &   & \vdots \\
\vdots & \vdots & \lambda_{s-1} &   & \vdots \\
0 & \lambda_{s-1} & \lambda_s - \alpha\, \lambda_{-1} & \ldots & \alpha\, \lambda_{s-2} \\
\alpha & \lambda_s & 2 \alpha\, \lambda_0 & \ldots & 2 \alpha\, \lambda_{s-1}
\end{array}\right].$$
Substracting a suitable multiple of the first column to each remaining column, we get:
$$ \det(M_\alpha) = \begin{array}{|ccccc|}
1 & 0 & 0 & \ldots & 0 \\
0 & \lambda_0 & \lambda_1 & \ldots & \lambda_s - \alpha\, \lambda_{-1} \\
\vdots & \vdots & \vdots & \ldots & \alpha\, \lambda_0 \\
\vdots & \vdots & \vdots &   & \vdots \\
\vdots & \vdots & \lambda_{s-1} &   & \vdots \\
0 & \lambda_{s-1} & \lambda_s - \alpha\, \lambda_{-1} & \ldots & \alpha\, \lambda_{s-2} \\
\alpha & \lambda_s  - \alpha\, \lambda_{-1} & \alpha\, \lambda_0 & \ldots & \alpha\, \lambda_{s-1}
\end{array}.$$
Then, using Laplace's expansion formula along the first row, it follows:
$$ \det(M_\alpha) = \begin{array}{|cccc|}
\lambda_0 & \lambda_1 & \ldots & \lambda_s - \alpha\, \lambda_{-1} \\
\vdots & \vdots & \ldots & \alpha\, \lambda_0 \\
\vdots & \vdots &   & \vdots \\
\vdots & \lambda_{s-1} &   & \vdots \\
\lambda_{s-1} & \lambda_s - \alpha\, \lambda_{-1} & \ldots & \alpha\, \lambda_{s-2} \\
\lambda_s  - \alpha\, \lambda_{-1} & \alpha\, \lambda_0 & \ldots & \alpha\, \lambda_{s-1}
\end{array}.$$

Let $k$ be the minimal integer such that $\lambda_k \neq 0$. Assume first that $k \neq -1$. The determinant of $M_\alpha$ takes the form:
$$ \det(M_\alpha) = 
\begin{array}{cc}
\underbrace{ \begin{array}{|cccc}
0 & \dots & 0 & \lambda_k \\
\vdots &  & \revddots & * \\
0 & \revddots &   & \vdots \\ 
\lambda_k & * & \ldots  & * \\ 
* & \ldots & \ldots & * \\
\vdots & & & \vdots \\
\vdots & & & \vdots \\ 
* & \ldots & \ldots & *  
\end{array} }_{k+1\ \mbox{\tiny columns}}
 & 
\underbrace{ \begin{array}{cccc|}
 * & \ldots & \ldots & * \\
\vdots & & & \vdots \\
\vdots & &  & \vdots \\
* & \ldots & \ldots & * \\
0 & \dots & 0 & \alpha\, \lambda_k \\
\vdots &  & \revddots & * \\
0 & \revddots &   & \vdots \\
\alpha\, \lambda_k & * & \ldots  & *
\end{array} }_{s-k\ \mbox{\tiny columns}}
\end{array} \ , $$
where ``$*$'' denote coefficients of the form $\alpha\, \lambda_i$ in the lower right diagonal block, and of the form $\lambda_i$ in all other three blocks. From this expression, it is quite standard to check that $\det(M_\alpha) = \pm (\alpha\, \lambda_k)^{s-k} \cdot \det(A) + \ldots$, where ``$\ldots$'' are polynomial expressions of degree $< s-k$ in $\alpha$ and $A$ is the upper left block. By assumption we have $\lambda_{k} \neq 0$, so that the component of degree $s-k$ of this expression is $\pm \lambda_k^{s+1} \alpha^{s-k} \neq 0$. Since $\FF$ is algebraically closed, there exists $\alpha \in \FF$ such that  $\det(M_\alpha) \neq 0$, and henceforth $\FF\langle x,y_\alpha \rangle = W(1,n)$.

In the case $k = -1$, a similar (and slightly easier) argument would show that $\det(M_\alpha) = (- \lambda_{-1})^{s+1} \alpha^{s+1} + \ldots$, leading to the same conclusion. This proves the proposition.

\subsubsection{\bf General case.} $ $ \\

\begin{thm}
Let $m \geq 1$ and $\nn \in \NN^m$ a multi-index with non-zero entries. The Lie algebra $W(m,\nn)$ is generated by one and a half elements if and only if $m = 1$. \\
\end{thm}

\proof\ The fact that $W(1,n)$ is generated by one and a half elements is Proposition \ref{zassenhaus}. Now assume that $m > 1$; we will prove that $W(m,\nn)$ is not generated by one and a half elements. By Proposition \ref{plongement}, there is an inclusion $W(m,\nn) \subseteq W(|\nn|,\un)$ such that $W(m,\nn)_s \subseteq W(|\nn|,\un)_{s'}$, where $W(m,\nn)_s$ and $W(|\nn|,\un)_{s'}$ are the components of highest degree of $W(m,\nn)$ and $W(|\nn|,\un)$ respectively. Let $x \in W(m,\nn)_s$ and $y \in W(m,\nn)$. Let $L = \FF \langle x,y \rangle \subseteq W(m,\nn)$ be the subalgebra generated by $x$ and $y$. Considering that $x \in W(|\nn|,\un)_{s'}$ and $y \in W(|\nn|,\un)$, we can apply Corollary \ref{cork} and get: 
$$\dim\, [L,L] \leq p^{|\nn|} < m p^{|\nn|} = \dim\, W(m,\nn).$$ 
Thus, $[L,L] \neq W(m,\nn)$, forcing also $L \neq W(m,\nn)$. The theorem is proved.

\section*{Acknowledgements}

This research was initiated during the author's stay at Bar-Ilan University in 2006. He would like to thank professor B. Kunyavski\u\i \ for drawing his attention to the question of one-and-a-half generation in Lie algebras, and for valuable discussions and comments.

The paper was written while the author was visiting the Collaborative Research Centre 701 at the University of Bielefeld. He would like to take the opportunity to thank Professor Claus Michael Ringel for his hospitality and support, and Professor Rolf Farnsteiner for fruitful exchanges during the redaction of this article.

He is also indebted to Professor Alexander Premet for explaining to him how to use the arguments of \ref{GeneratingSetsGradedCase} to deal with the non-graded simple Lie algebras. Finally, he would like to thank the referee for careful reading and very valuable comments.

\end{document}